\documentclass[reqno]{amsart}

\usepackage{geometry}
\usepackage[english]{babel}
\usepackage{latexsym, mathrsfs, verbatim}
\usepackage{amsfonts, amsthm, amsmath, amssymb, bm}
\usepackage{hyperref}
\usepackage{amsaddr}

\newtheorem{prop}{Proposition}[section]
\newtheorem{lem}[prop]{Lemma}
\newtheorem{cor}[prop]{Corollary}
\newtheorem{thm}[prop]{Theorem}

\theoremstyle{definition}
\newtheorem{rem}[prop]{Remark}
\newtheorem{defi}[prop]{Definition}
\newtheorem{ex}[prop]{Example}

\def\Z{\mathbb{Z}}

\def\equad{\quad \textrm{ and } \quad}

\def\H{\mathrm{H}}
\def\SMA{\mathrm{SMA}}
\def\SMAS{\mathrm{SMAS}}
\def\IHS{\mathrm{IHS}}

\def\Mo{{}^0\mathrm{MPF}}

\def\T{\mathsf{t}}

\def\E{\mathcal{E}}

\def\diag{\mathrm{diag}}
\def\dcup{\sqcup}

\numberwithin{equation}{section}

\usepackage[table,dvipsnames]{xcolor}
\usepackage[normalem]{ulem}
 \usepackage{datetime}
 
 \def\sha{\cellcolor{black!5}}
 \def\shb{\cellcolor{black!10}}
\def\shc{\cellcolor{black!20}}
\def\shd{\cellcolor{black!25}}

\hypersetup{ 
 colorlinks =true,
 linkcolor=Sepia,
 citecolor=Plum,
 urlcolor=olive,
 }

\begin{document}

\title{Signed magic arrays: existence and constructions}

\author[Fiorenza Morini]{Fiorenza Morini}
\address{Dipartimento di Scienze Matematiche, Fisiche e Informatiche, Universit\`a di Parma,\\
Parco Area delle Scienze 53/A, 43124 Parma, Italy}
\email{fiorenza.morini@unipr.it}

\author[Marco Antonio Pellegrini]{Marco Antonio Pellegrini}
\address{Dipartimento di Matematica e Fisica, Universit\`a Cattolica del Sacro Cuore,\\
Via della Garzetta 48, 25133 Brescia, Italy}
\email{marcoantonio.pellegrini@unicatt.it}

\begin{abstract}
Let $m,n,s,k$ be four integers such that $1\leqslant s \leqslant n$, $1\leqslant k\leqslant m$ and $ms=nk$. A signed magic array $\SMA(m,n; s,k)$ is an $m\times n$ partially filled  array
whose entries belong to the subset $\Omega\subset \Z$, where $\Omega=\{0,\pm 1, \pm 2,\ldots, \pm (nk-1)/2\}$ if $nk$ is odd and
$\Omega=\{\pm 1, \pm 2, \ldots, \pm nk/2\}$ if $nk$ is even, satisfying the following requirements:
$(a)$ every $\omega \in \Omega$ appears  once in the array;
$(b)$ each row contains exactly $s$ filled cells and each column contains exactly $k$ filled cells;
$(c)$ the sum of the elements in each row and in each column is $0$.
In this paper we construct these arrays when $n$ is even and $s,k\geqslant 5$ are coprime integers. 
This allows us to provide a complete answer to a problem posed in 2017 by Khodkar, Schulz and Wagner, giving the necessary and sufficient conditions for the existence of an $\SMA(m,n; s,k)$ for all admissible values of $m,n,s,k$.
\end{abstract}

\maketitle

\section{Introduction}

In recent years, many authors have worked with magic squares and various generalizations of these objects.
For instance, one can consider magic  rectangles \cite{H1,H2}, magic  hypercubes \cite{S}, or magic squares on abelian groups
\cite{WS}.
In particular, Abdollah Khodkar, Christian Schulz and Nathan Wagner introduced in \cite{KCW} a class of partially filled arrays, called
\emph{signed magic arrays}.

\begin{defi}\label{def:Magic}
Let $m,n,s,k$ be four integers such that $1\leqslant s \leqslant n$, $1\leqslant k\leqslant m$ and $ms=nk$.
A \emph{signed magic array} $\SMA(m,n; s,k)$ is an $m\times n$ partially filled  array
with entries in $\Omega\subset \Z$, where $\Omega=\{0,\pm 1, \pm 2,\ldots, \pm (nk-1)/2\}$ if $nk$ is odd and
$\Omega=\{\pm 1, \pm 2, \ldots, \pm nk/2\}$ if $nk$ is even, such that
\begin{itemize}
\item[{\rm (a)}] every $\omega \in \Omega$ appears  once in the array;
\item[{\rm (b)}] each row contains exactly $s$ filled cells and each column contains exactly $k$ filled cells;
\item[{\rm (c)}] the sum of the elements in each row and in each column is $0$.
\end{itemize}
\end{defi}

Replacing in the previous definition the subset $\Omega$ with any subset $\Psi\subset\Z$ such that $\{\Psi,-\Psi\}$ is a partition of
$\{\pm 1, \pm 2,\ldots,\pm nk\}$, one obtains an \emph{integer} Heffter array $\H(m,n;s,k)$, an object introduced by Dan Archdeacon in \cite{A}.
In fact, as shown in \cite{MP4}, magic rectangles, signed magic arrays and Heffter arrays are all members of the broader family of the \emph{magic partially filled arrays}.
Speaking of integer Heffter arrays, we recall 
that Archdeacon conjectured  (see \cite[Conjecture 6.3]{A}) that there 
exists an integer Heffter array $\H(m, n; s, k)$ for all $m, n, s, k$ with $s, k \geqslant 3$
and $ms = nk$, provided that the necessary condition $ms\equiv 0,3 \pmod 4$ is satisfied.
The reader interested in constructing these objects can find a survey on almost all known results in \cite{DP}. 
Further methods for constructing rectangular integer Heffter arrays have been recently described in \cite{MP3,PT}.
These methods, adapted here to work with signed magic arrays, allow us to completely solve the existence problem of an $\SMA(m,n;s,k)$.

\begin{thm}\label{main}
Let  $m,n,s,k$ be four integers such that $1\leqslant s\leqslant n$, $1\leqslant k \leqslant m$ and $ms=nk$. 
There exists an $\SMA(m,n;s,k)$ if and only if 
one of the following cases occurs:
\begin{itemize}
\item[(1)] $k=m=s=n=1$;
\item[(2)] $k=2$, $m=2$ and $s=n\equiv 0,3 \pmod 4$;
\item[(3)] $k=2$ and $m,s\geqslant 3$;
\item[(4)] $s=2$, $n=2$ and $k=m\equiv 0,3 \pmod 4$;
\item[(5)] $s=2$ and $n,k\geqslant 3$;
\item[(6)] $s,k\geqslant 3$.
\end{itemize}
\end{thm}

Clearly, if $A$ is an $\SMA(m,n;s,k)$, then its transpose $A^\T$ is an $\SMA(n,m;k,s)$.
Furthermore, an $\SMA(m,n;s,1)$ exists only when $m=s=n=1$, in which case the matrix one is looking for is
$\begin{array}{|c|}\hline 0 \\ \hline\end{array}$.
The existence problem of an $\SMA(m,n;s,2)$ and of an $\SMA(m,n;s,3)$
has been already solved in \cite{KE} and in \cite{KLE}, respectively.
Constructions of signed magic arrays 
$\SMA(m,n;s,k)$ have been described firstly assuming that $s,k\geqslant 4$ are even integers (see \cite[Theorem 1.6]{MP2}), and then assuming the weaker hypothesis that  $\gcd ( s , k ) \geqslant 2$
(see \cite[Theorem 1.9]{MP3}). The case when 
$s \equiv  0 \pmod 4$ has been dealt with  in \cite[Theorem 1.10]{MP3}, while 
the case when $nk$ is odd has been completely solved in \cite[Corollary 3.5]{MP4}.
We can resume the results obtained in \cite{MP2,MP3, MP4} as follows.

\begin{prop}\label{fatto}
Let  $m,n,s,k$ be four integers such that $3\leqslant s\leqslant n$, $3\leqslant k \leqslant m$ and $ms=nk$.
There  exists an $\SMA(m,n;s,k)$ in each of the following cases:
\begin{itemize}
\item[(1)] $\gcd(s,k)\geqslant 2$;
\item[(2)] $s\equiv 0 \pmod 4$;
\item[(3)] $nk$ is odd.
\end{itemize}
\end{prop}

Once more, we remark that all the previous results on signed magic arrays \cite{KE,KLE,MP2,MP3,MP4} have been obtained with constructive methods.
So, to prove Theorem \ref{main} it suffices to consider the case when $nk$ is even, $s,k\geqslant 5$ and $\gcd(s,k)=1$.
Under these assumptions, $s$ and $k$ cannot be two even integers, so we may assume that $k$ is odd
and $s\not \equiv 0\pmod 4$.
If $s$ and $k$ are  odd coprime integers, from $ms=nk$ even, it follows that $m=2ck$ and $n=2cs$ for some
positive integer $c$.
On the other hand, if $s\equiv 2 \pmod 4$, we can write $s=2\bar s$ for some odd integer $\bar s\geqslant 3$ such that
$\gcd(\bar s, k)=1$. Again, from $ms=nk$, we obtain $m=ck$ and $n=2c\bar s$ for some positive integer $c$.
In both cases, we prove the existence of an $\SMA(m,n;s,k)$ by constructing a signed magic array set, see \cite{MP4}.

\begin{defi}\label{SMAS}
A \emph{signed magic array set} $\SMAS(a,b;  e)$ is a set of
$e$ arrays of size $a\times b$ with entries in $\Omega\subset \Z$,
where $\Omega=\{0,\pm 1, \pm 2,\ldots, \pm (abe-1)/2\}$ if $abe$ is odd and
$\Omega=\{\pm 1, \pm 2, \ldots, \pm abe/2\}$ if $abe$ is even,
such that
\begin{itemize}
\item[{\rm (a)}] every $\omega \in \Omega$ appears  once and in a unique array;
\item[{\rm (b)}] for every array, the sum of the elements in each row and in each column is $0$.
\end{itemize}
\end{defi}

\begin{figure}[ht]
\begin{footnotesize}
$$\begin{array}{|c|c|c|c|c|c|c|}\hline
\sha 10 & \sha -4  &\sha  -6   & \shc 5   & \shc  -5  &  \shd  7 & \shd -7 \\\hline
\sha -2 & \sha 19  & \sha -17  & \shc  -32 & \shc  32  & \shd -33 &\shd  33 \\\hline
\sha -8 & \sha -15 &\sha  23   & \shc  27  & \shc  -27 & \shd 26 &\shd  -26 \\\hline
\shb 1 &\shb -30  &\shb 29   & 14  & -16 & -18 & 20 \\\hline
\shb -1 &\shb 30  &\shb -29  & -14 & 16  & 18 & -20 \\\hline
\end{array},\quad \begin{array}{|c|c|c|c|c|c|c|}\hline
\sha -10 & \sha 4    & \sha 6  &  \shc 9   & \shc  -9  &\shd  11  &\shd  -11 \\\hline
\sha 2   & \sha -19 & \sha 17  &  \shc -34 & \shc  34  &\shd  -35 &\shd  35 \\\hline
\sha 8   & \sha 15  & \sha -23 &  \shc 25  & \shc  -25 &\shd  24  &\shd  -24 \\\hline
\shb 3   & \shb-31 &\shb 28  & 12  & -13 & -21 & 22 \\\hline
\shb -3  &\shb 31  &\shb -28 & -12 & 13  & 21  & -22 \\\hline
\end{array}.$$
\end{footnotesize}

\caption{An $\SMAS(5,7;2)$.}\label{572}
\end{figure}

In the next sections, we prove the following results that, together with Proposition \ref{fatto}, will imply the validity of Theorem \ref{main}, as explained in Section \ref{concl}.

\begin{prop}\label{p:sets}
Given a positive integer $c$ and two odd integers $a,b\geqslant 5$ such that $(a,b)\neq (5,5)$,
there  exists an $\SMAS(a,b; 2c)$.
\end{prop}

\begin{prop}\label{6}
Given an odd integer $b\geqslant 5$,
there exists an $\SMAS(6,b;c)$ for all  $c\geqslant 1$.
\end{prop}

We recall that, when $nk$ is even, an $\SMA(m,n;s,k)$ is an integer $2$-fold Heffter array  ${}^2\H(m,n;s,k)$  with some additional properties.

\begin{defi}\cite{CP,MP2}
Let $m,n,s,k$ be four integers such that $1\leqslant  s\leqslant n$, $1\leqslant k \leqslant m$ and
$ms = nk$. 
Set $v=nk+1$ and $\Phi=\left\{1, \ldots,\left\lfloor \frac{v}{2} \right\rfloor \right\}\subset \Z$.
An \emph{integer $2$-fold Heffter array} ${}^2\H(m,n;s,k)$ is an $m\times n$ partially filled array with elements from $\Phi \cup -\Phi$ such that
\begin{itemize}
\item[{\rm (a)}] each row contains exactly $s$ filled cells and each column contains exactly $k$ filled cells;
\item[{\rm (b)}] denoting by $\rho(x)$ the number of the entries whose absolute value is equal to a positive integer $x$, if $v$ is odd, then $\rho( x ) = 2$ for every $x\in\Phi$; if $v$ is even, then
$\rho(x)=2$ for every $x \in \Phi \setminus \left\{ \frac{v}{2}  \right\}$ while $\rho\left(\frac{v}{2}\right)
=1$;
\item[{\rm (c)}] the sum of the elements in each row and in each column is $0$.
\end{itemize}
\end{defi}

As a consequence of Theorem \ref{main}, we have the following.

\begin{cor}\label{fold}
Let $m,n,s,k$ be four integers such that $2\leqslant  s\leqslant n$, $2\leqslant k \leqslant m$ and
$ms = nk$. Suppose $nk\not \equiv 3 \pmod 4$.
There exists an integer ${}^2\H(m,n;s,k)$ if and only if $nk$ is even and
$$(m,n,s,k)\not \in \left\{ (\ell, 2,2,\ell), (2,\ell,\ell, 2): \ell\geqslant 1 \text{ is such that } \ell \equiv 1,2 \pmod 4\right\}.$$
\end{cor}

Signed magic arrays  can be interpreted as zero-sum magic partially filled arrays, according to the following.

\begin{defi}\cite{MP4}
A \emph{zero-sum magic partially filled array} $\Mo_\Omega(m,n; s,k)$ on a subset $\Omega$ of an abelian group $(\Gamma,+)$
is a  partially filled array of size $m\times n$ with entries in $\Omega$ such that
\begin{itemize}
\item[{\rm (a)}] every $\omega \in \Omega$ appears once in the array;
\item[{\rm (b)}] each row contains $s$ filled cells and each column contains $k$ filled cells;
\item[{\rm (c)}] the sum of the elements in each row  and in each column is  $0_\Gamma$.
\end{itemize}
\end{defi}

The entries of  an $\SMA(m,n;s,k)$   can be also viewed as elements of a finite cyclic group. Setting $\Gamma=\Z_{nk}$ if $nk$ is odd and $\Gamma=\Z_{nk+1}$ if $nk$ is even,  an $\SMA(m,n;s,k)$  is a 
zero-sum magic partially filled array on $\Gamma$ if $nk$ is odd and on 
$\Gamma\setminus\{0_\Gamma\}$ if $nk$ is even.
Consequently, signed magic arrays may prove to be potentially useful for graph theorists working on graph labeling. We now briefly describe a possible connection between these arrays (and signed magic array sets) and magic labelings (see \cite[Sections 5.1 and 5.6]{G}). 

A bipartite biregular graph $G(V,E)$ is a graph whose vertex set can be written as a disjoint union $V=V_1\cup V_2$ with
$|V_1|=m$, $|V_2|=n$, and where each vertex of $V_1$ is connected with exactly $s$ vertices of $V_2$,
and each vertex of $V_2$ is connected with exactly $k$ vertices of $V_1$.
Now, let $M$ be an $\SMA(m,n;s,k)$.
We associate to  $M$ a bipartite biregular graph $\Phi_M=G(V,E)$ by taking a set $V_1$ of $m$ points, a set $V_2$ of $n$ points, and drawing an edge
$e_{i,j}$ between the $i$-th vertex of $V_1$ and the $j$-vertex of $V_2$ if the cell $(i,j)$ of $M$ is not empty.
Define the labeling $f_M: E\to \Gamma$, where 
$f_M(e_{i,j})$ is the entry of the corresponding cell $(i,j)$ of~$M$.

A \emph{$\Z_{nk}$-supermagic labeling} of a graph $G(V, E)$ with $|E| = nk$ is a bijection from $E$ to
$\Z_{nk}$ such that the sum of labels of all incident edges of every vertex
$v \in V$ is equal to the same element $x \in \Z_{nk}$, see \cite{FMMM}.
If $M$ is an $\SMA(m,n;s,k)$ with $nk$ odd, then function $f_M$ is a $\mathbb{Z}_{nk}$-supermagic labeling of $\Phi_M$, where the constant $x$ is $0_{\Z_{nk}}$. 

A graph $G$ is said to be \emph{zero-sum $\Z_{nk+1}$-magic} if there exists a labeling of the edges of $G$ with
elements of $\Z_{nk+1}\setminus \{0\}$ such that, for each vertex $v$, the sum of the labels of the edges incident with $v$ is equal to $0_{\Z_{nk+1}}$, see \cite{ARZ}.
If $M$ is an $\SMA(m,n;s,k)$ with $nk$ even, then $\Phi_M$ is a zero-sum
$\Z_{nk+1}$-magic graph.

Finally, given a graph $G=(V,E)$, a  bijection  $f: V \to \Z_ {abe+1} \setminus \{0_{\Z_{abe+1}}\}$  is said to be a
\emph{$\Z_{abe+1}$-distance magic labeling} of $G$ if there exists $\mu\in \Z_{abe+1}$
such that $\omega(v) = \sum\limits_{u \in N(v)} f(u)$ is $\mu$ for all vertices $v\in V$, see \cite{C1,C2}.
If $S$ is an $\SMAS(a,b;e)$ with $abe$ even, then we can consider the complete $a$-partite graph $K_{be\times a}$ consisting of $be$ parts of the same cardinality $a$. Labeling the vertices of each part with the $a$ entries of the $be$ columns of the elements of $S$, we obtain a $\Z_{abe+1}$-distance magic labeling of $K_{be\times a}$ where the constant $\mu$ is $0_{\Z_{abe+1}}$.
\smallskip

Explicit constructions of our signed magic array sets are given in Sections~\ref{gen}, \ref{sm} and \ref{s6}. The reader interested only in the proofs of our main results can refer directly to Section~\ref{concl}.
Those, who want to understand how we have truly solved this decade-long problem, will find similar constructions (and statements) each requiring small adjustments or different approaches, due to the presence of various parameters.  We have  unified some of them while trying to maintain a certain level of readability. To help these readers,
explicit algorithmic implementations in \textsc{gap} are available up to request to the authors.

\section{Notation and preliminary results}\label{basic}

Given an array $A$, we denote by 
$\E(A)$ the list of its entries;
$\sigma_r(A)$ and $\sigma_c(A)$ denote, respectively, the sequence of the sums of the elements of each row and column of $A$.
If $\sigma_r(A)=(0,\ldots,0)$ and $\sigma_c(A)=(0,\ldots,0)$,
the array $A$ is said to be a zero-sum block.
Given a set $\mathfrak{S} = \{ A_1, A_2 ,
\ldots, A_r\}$ of arrays, we set 
$\E ( \mathfrak{S} ) = \cup_i \E (A_i )$.

Given an integer $d\geqslant 1$, if $a,b$ are two integers such that $a\equiv b \pmod{d}$, then we use the notation
$$[a,b]_d =\left\{a+id\mid 0\leqslant i \leqslant \frac{b-a}{d}\right\}\subset \Z,$$ whenever 
$a\leqslant b$. If $a>b$, then $[a,b]_d=\varnothing$.
If $d=1$, we simply write $[a,b]$.
Given a subset $\Omega$ consisting of positive integers, 
we write $\pm \Omega=\{\pm \omega: \omega \in \Omega \}$.
Given two positive integers $\varepsilon$ and $x$,
a set $\{x,x+\varepsilon\}$ will be called a $2$-set of type $\varepsilon$.

Following \cite{PT}, we consider the following objects.

\begin{defi}\label{IHS}
An \emph{integer Heffter array set} $\IHS(a,b; c)$ is a set of $c$ arrays of size $a\times b$  such that
\begin{itemize}
\item[{\rm (a)}] the entries of all these arrays constitute a subset $\Omega\subset \Z$  such that $\{\Omega, -\Omega\}$ is a partition of 
$\{\pm 1, \pm2, \ldots, \pm abc \}$;
\item[{\rm (b)}] every $\omega \in \Omega$ appears  once and in a unique array;
\item[{\rm (c)}] for every array, the sum of the elements in each row  and in each column is $0$.
\end{itemize}
\end{defi}

\begin{thm}\label{thIHS}\cite{PT}
Let $a,b,c$ be positive integers such that
$abc\equiv 0,3 \pmod 4$.
Suppose that $a,b\geqslant 7$ are odd integers.
Then,  there exists an $\IHS(a,b;c)$.
\end{thm}

If $\{A_1,A_2,\ldots,A_c\}$ is an  $\IHS(a,b;c)$,
then $\{A_1,-A_1,A_2,-A_2,\ldots,A_c,-A_c\}$
is clearly an $\SMAS(a,b;2c)$.
Hence, from Theorem \ref{thIHS} we obtain the following.

\begin{cor}\label{IHS->SMAS}
Let $a,b,c$ be positive integers such that
$abc\equiv 0,3 \pmod 4$.
Suppose that $a,b\geqslant 7$ are odd integers.
Then,  there exists an $\SMAS(a,b;2c)$. 
\end{cor}

\begin{rem}\label{tr}
An $\SMAS(a,b;e)$ exists if and only if there exists an $\SMAS(b,a;e)$.
\end{rem}

In the next two sections, we will construct 
$\SMAS(a,b;2c)$ when either $a=5$ or $abc \equiv 1,2 \pmod 4$. 
Roughly speaking, to build the elements of our signed magic array set we will use some basic zero-sum blocks of size $3\times 3$, $2\times 3$, $2\times 4$ or $2\times 6$ 
(see, for instance, Figure \ref{572}). To this purpose, we start with the following constructions.

\begin{lem}\label{BlockA}
Given three nonnegative integers $\alpha,\beta,u$ such that
$u\leqslant \alpha+1 $,
there exists a set $\mathfrak{A}=\mathfrak{A}(\alpha,\beta,u)$ consisting of $2u$ square zero-sum blocks of size $3$ such that
$\E(\mathfrak{A})$ consists of the following disjoint subsets:
$$\begin{array}{rcl}
\E(\mathfrak{A})  & = & \pm[2\alpha+2, 2\alpha+8u]_2 \dcup 
\pm [4\alpha+6u+4,4\alpha+10u]_4 \dcup\\
&& \pm [2\beta+1, 2\beta+4u-1]_2\dcup  \pm [2\alpha+2\beta+2u+3,  2\alpha+2\beta+10u-1]_4.
\end{array}$$
\end{lem}

\begin{proof}
For all $i\in [0,u-1]$, define
$$A_{i} =  \begin{array}{|c|c|c|}\hline
 4\alpha+6u+4+4i     &  -(2\alpha+2u+2+2i) &  -(2\alpha+4u+2+2i) \\ \hline
-(2\alpha+2+2i)   & 2\alpha+2\beta+2u+3+4i  &  -(2\beta+2u+1+2i) \\\hline
-(2\alpha+6u+2+2i) & -(2\beta+1+2i) & 2\alpha+2\beta+6u+3+4i\\ \hline
\end{array}.$$
The set $\mathfrak{A}=\{+A_{i},-A_{i} :  i \in [0,u-1]\}$ has the required properties.
\end{proof}

\begin{lem}\label{BlockB}
Given three positive integers $\alpha,\beta,u$ such that
$\beta\geqslant \alpha+3u-2$,
there exists a set $\mathfrak{B}=\mathfrak{B}(\alpha,\beta,u)$ consisting of $u$ zero-sum blocks of size $2\times 3$ such that
$\E(\mathfrak{B})$ consists of the following disjoint subsets:
$$\E(\mathfrak{B})  = \pm [\alpha,\alpha+2u-2]_2\dcup
\pm [\beta-u+1,\beta] \dcup
\pm [\alpha+\beta,\alpha+\beta+u-1].$$
\end{lem}

\begin{proof}
For every $j\in [0, u-1]$, define
$$B_j=\begin{array}{|c|c|c|}\hline
\alpha +2j & -( \alpha+\beta+j) & \beta-j \\ \hline
-(\alpha+2j) & \alpha+\beta+j & -(\beta-j ) \\\hline
\end{array}.$$
The set $\mathfrak{B}=\left\{B_j: j \in [0, u -1]\right\}$ has the required properties.
\end{proof}

Let $Q_1=Q_1(\{x_1,x_1+1\},\{x_2,x_2+1\})$, $Q_2=Q_2(\{y_1,y_1+2\},\{y_2,y_2+2\})$,
$Q_3=Q_3(\{z_1,z_1+4\},\{z_2,z_2+4\})$, $R_1=R_1(\{ y_1, y_1+2\}, \{x_1,x_1+1\},\{x_2,x_2+1\})$ 
and $R_2=R_2(\{ z_1, z_1+4\}, \{y_1,y_1+2\},\{y_2,y_2+2\})$ be the following zero-sum blocks:
$$
\begin{array}{rcl}
Q_1 & =& \begin{array}{|c|c|c|c|}\hline
 x_1 & -(x_1+1) & -x_2 & x_2+1 \\ \hline 
 -x_1 & x_1+1 & x_2 & -(x_2+1) \\ \hline
\end{array},\\\\[-8pt]
Q_2 & =& \begin{array}{|c|c|c|c|}\hline
 y_1 & -(y_1+2) & -y_2 & y_2+2 \\ \hline 
 -y_1 & y_1+2 & y_2 & -(y_2+2) \\\hline
\end{array}, \\\\[-8pt]
Q_3 & =& \begin{array}{|c|c|c|c|}\hline
 z_1 & -(z_1+4) & -z_2 & z_2+4 \\ \hline 
 -z_1 & z_1+4 & z_2 & -(z_2+4) \\\hline
\end{array}, \\\\[-8pt]
R_1 & =& \begin{array}{|c|c|c|c|c|c|}\hline
y_1 & -(y_1+2) & -x_1 & x_1+1 & -x_2 & x_2+1 \\ \hline 
-y_1 & y_1+2 & x_1 & -(x_1+1) & x_2 & -(x_2+1) \\\hline
\end{array},\\\\[-8pt]
R_2 & =& \begin{array}{|c|c|c|c|c|c|}\hline
z_1 & -(z_1+4) & -y_1 & y_1+2 & -y_2 & y_2+2 \\ \hline 
-z_1 & z_1+4 & y_1 & -(y_1+2) & y_2 & -(y_2+2) \\\hline
\end{array}.
\end{array}$$

\begin{lem}\label{BlockC}
Given three positive integers $\alpha,\beta,\gamma$ such that
$\beta >\gamma+3 > \alpha+12 > 22$, there exists 
a set $\mathfrak{C}=\mathfrak{C}(\alpha,\beta,\gamma)$ consisting of two square blocks of size $5$ such that
 $$\sigma_r(C)=( 0,0,0,1,-1) \equad \sigma_c(C)=(0,0,0,0,0) \text{ for all } C \in \mathfrak{C}$$
 and $\E(\mathfrak{C})$ consists of the following disjoint subsets:
 $$\E(\mathfrak{C})   =\pm [1,7]_2 \dcup \pm [2,10]_2 
 \dcup  
\pm \{\alpha, \alpha+4,\alpha+6, \alpha+9\} \dcup 
 \pm [\gamma, \gamma+3]\dcup 
\pm [\beta, \beta+7]. $$
\end{lem}

\begin{proof}
Let $C_1$ and $C_2$ be the following arrays: 
$$\begin{array}{|c|c|c|c|c|}\hline
10 &      -8    &          -2  &  5         &    -5   \\ \hline
-6 & -\alpha   & \alpha+6      & -(\beta+6) & \beta+6 \\ \hline
-4 & \alpha+9  &  -(\alpha+4)  & \beta      & -(\beta+1) \\ \hline
 1 & -(\beta+5)  &   \beta+3   & -(\gamma+1) & \gamma+3 \\ \hline
-1 &   \beta+4   & -(\beta+3)  & \gamma+2 & -(\gamma+3) \\ \hline
  \end{array},$$
$$\begin{array}{|c|c|c|c|c|}\hline
-10 &      8    &          2   & 7          &    -7 \\ \hline 
6 & \alpha   & -(\alpha+6)     & -(\beta+7) &  \beta+7  \\ \hline 
4 & -(\alpha+9)  &  \alpha+4   & \beta+1    & -\beta   \\ \hline 
-3 & \beta+5  &   -(\beta+2)   & \gamma+1   & -\gamma \\ \hline
 3 & -(\beta+4)  & \beta+2     & -(\gamma+2) & \gamma \\ \hline
  \end{array}.$$
The set $\mathfrak{C}=\{C_1,C_2\}$ has the required properties.
  \end{proof}

The blocks in the previous set $\mathfrak{C}$ 
will require the use of  some
special $2\times 4$ and $2\times 2$ blocks.
Let $U_1(\{y_1,y_1+2\}, \{x_1,x_1+1\})$ and $U_2(\{x_1,x_1+1\})$ be the blocks
$$U_1=\begin{array}{|c|c|c|c|}\hline
y_1 & - (y_1+2) & -x_1 & x_1+1 \\ \hline 
-y_1 & y_1+2 &  x_1 & -(x_1+1)\\ \hline
\end{array} \equad 
U_2=\begin{array}{|c|c|}\hline
 x_1 & -(x_1+1)\\ \hline
-x_1 & x_1+1 \\ \hline 
\end{array}.$$
We have $\sigma_r(U_1)=\sigma_r(U_2)=(-1, +1)$, $\sigma_c(U_1)=(0,0,0,0)$ and $\sigma_c(U_2)=(0,0)$.
\smallskip

\section{The general case}\label{gen}

In this section we construct an $\SMAS(a,b;2c)$ when  $a,b\geqslant 5$ are odd integers such that $(a,b)\neq (5,5),(5,7),(7,5)$, and $c$ is a positive integer such that $c\not \equiv 0 \pmod 4$.

\subsection{Subcase \texorpdfstring{$c\equiv 3 \pmod 4$}{c ≡ 3 (mod 4)}}

\begin{lem}\label{59-3}
Suppose that $b\geqslant 9$ is such that $b\equiv 1 \pmod 4$.
There exists an $\SMAS(5,b;2c)$ for all $c\geqslant 3$ such that $c\equiv 3 \pmod 4$.
\end{lem}

\begin{proof}
Write $c=4t+3$ where $t\geqslant 0$.
Our $\SMAS(5, b;2c)$ will consist
of $8t+6$ arrays of type
\begin{equation}\label{59}
\begin{array}{|c|c|c|c|c|c|c|c|c|}\hline
  X_3      &   X_2^\T & X_{2}^\T  & X_{2}^\T 
  &  X_2^\T & X_{2}^\T  & \cdots & X_2^\T & X_{2}^\T \\ \hline
  X_{2}    &  \multicolumn{3}{c|}{Z_6} 
  &  \multicolumn{2}{c|}{Z_4} & \cdots &  \multicolumn{2}{c|}{Z_4} 
   \\ \hline  
\end{array},
\end{equation}
where each $X_r$ is a zero-sum block of size $r \times 3$ and each $Z_{\ell}$ is a zero-sum block of size $2\times \ell$.

Write $b=4w+9$, where  $w\geqslant 0$.
Hence, we will construct
$8t+6$ blocks $X_3$, $4(w+2)(4t+3)$
blocks $X_2$, $2w(4t+3)$ blocks $Z_4$ and $8t+6$ blocks  $Z_6$.
So, take
$$\begin{array}{rcl}
\mathfrak{A}  &=& \mathfrak{A}( 4t+2,\;  4(w+2)(4t+3) ,\; 4t+3 ),\\
\mathfrak{B} & =& 
\mathfrak{B}(1,\; 4(16w+37)t+ 48w+111  , \; 4(w+2)(4t+3) )
\end{array}$$
as in Lemmas \ref{BlockA} and \ref{BlockB}, respectively.
Then, 
$$\E(\mathfrak{A})\dcup \E(\mathfrak{B})=\pm [1, 20w(4t+3)+180t+135]\setminus \pm \left(\mathcal{M}_1\dcup \ldots \dcup \mathcal{M}_6 \right),$$
where
$$\begin{array}{rcl}
\mathcal{M}_1 & =&  [2,8t+4]_2, \\
\mathcal{M}_2 & =& [40t+32, 56t+36]_4,\\
\mathcal{M}_3 & =& [56t+40,  8w(4t+3)+112t+78]_2,\\
\mathcal{M}_4 & =& [8w(4t+3)+80t+63,   8w(4t+3)+112t+75]_4,\\
\mathcal{M}_5 & =& [8w(4t+3)+112t+79, 8w(4t+3)+112t+80],\\
\mathcal{M}_6 & =& [8w(4t+3)+112t+82, 12w(4t+3)+116t+87].
\end{array}$$
So, we replace the $8t+6$ instances of $X_3$ with the arrays in  $\mathfrak{A}$ and the $4(w+2)(4t+3)$ instances of $X_2$ with the arrays in $\mathfrak{B}$.

Next, write
$$\mathcal{M}_1=\mathcal{M}_1^1\dcup \mathcal{M}_1^2\dcup \mathcal{M}_1^3
\equad
\mathcal{M}_6=\mathcal{M}_6^1\dcup \mathcal{M}_6^2\dcup \mathcal{M}_6^3
\dcup\mathcal{M}_6^4\dcup \mathcal{M}_6^5, $$
where
$$\begin{array}{rcl}
\mathcal{M}_1^1 & =& [2,8t-2]_4,\\
\mathcal{M}_1^2 & =& [4,8t]_4,\\
\mathcal{M}_1^3 & =& [8t+2,8t+4]_2,\\
\mathcal{M}_6^1 & =& [8w(4t+3)+112t+82, 8w(4t+3)+112t+86  ]_4, \\
\mathcal{M}_6^2 & =& [8w(4t+3)+112t+83, 8w(4t+3)+112t+87 ]_4,\\
\mathcal{M}_6^3 & =& [8w(4t+3)+112t+84, 8w(4t+3)+112t+85],  \\
\mathcal{M}_6^4 & =& [8w(4t+3)+112t+88, 12w(4t+3)+ 116t+86]_2,\\
\mathcal{M}_6^5 & =& [8w(4t+3)+112t+89, 12w(4t+3)+ 116t+87]_2.
\end{array}$$

The set $\mathcal{M}_1^1\dcup \mathcal{M}_1^2 \dcup \mathcal{M}_2\dcup \mathcal{M}_4\dcup\mathcal{M}_6^1\dcup \mathcal{M}_6^2$  can be written as a disjoint union 
$F_1\dcup \ldots \dcup F_{\alpha}$, where $\alpha= 8t+5$ and each $F_i$ is
a $2$-set of type $4$.
The set $\mathcal{M}_1^3\dcup \mathcal{M}_3 \dcup \mathcal{M}_6^4\dcup \mathcal{M}_6^5$
can be written as a disjoint union 
$G_1\dcup \ldots \dcup G_{\beta}$, where $\beta=4w(4t+3)+16t+11$  and each $G_j$ is
a $2$-set of type $2$.
Finally, the set $\mathcal{M}_5\dcup \mathcal{M}_6^3$
can be written as a disjoint union 
$H_1\dcup H_2$, where each $H_h$ is a $2$-set of type $1$. 
Call
$$\begin{array}{rcl}
\mathfrak{R}_1 & =& \left\{R_2(F_i,G_{2i-1},G_{2i} ): i \in [1,8t+5] \right\},\\
\mathfrak{R}_2 & =& \left\{R_1(G_{16t+11},H_1,H_2 )\right\}, \\
\mathfrak{S}  & =& \left\{Q_2 (G_{2j},G_{2j+1}) : j \in [ 8t+6, 2w(4t+3)+8t+5] \right\} .
\end{array}$$
Then $\mathfrak{R}=\mathfrak{R}_1\dcup \mathfrak{R}_2$ consists of $8t+6$ zero-sum blocks of size $2\times 6$ and
$\mathfrak{S}$ consists of $2w(4t+3)$ zero-sum blocks of size $2\times 4$.
Furthermore, $\E(\mathfrak{R})\dcup \E(\mathfrak{S})=\pm 
(\mathcal{M}_1\dcup \ldots\dcup \mathcal{M}_6)$.
Hence, we replace the instances of $Z_6$ with the arrays in  $\mathfrak{R}$
and the instances of $Z_4$ with the arrays in~$\mathfrak{S}$.
\end{proof}

\begin{ex}
 We construct an $\SMAS(5,13;6)$ following the proof of the previous lemma.
Setting $t=0$ and $w=1$, we need
 $6$ blocks $X_3$, $36$ blocks $X_2$, $6$ blocks $Z_4$ and $6$ blocks $Z_6$. In particular, we have 
 $\mathcal{M}_1^1=\mathcal{M}_1^2=\varnothing$, $\mathcal{M}_2\dcup \mathcal{M}_4\dcup \mathcal{M}_6^1 \dcup \mathcal{M}_6^2
=[32,36]_4 \dcup [87,99]_4 \dcup [106,110]_4 \dcup [107,111]_4$, $\mathcal{M}_1^3\dcup\mathcal{M}_3\dcup \mathcal{M}_6^4
\dcup \mathcal{M}_6^5
=[2,4]_2 \dcup [40,102]_2 \dcup 
[112,122]_2 \dcup [113,123]_2$ and $\mathcal{M}_5\dcup \mathcal{M}_6^3=[103,104]\dcup [108,109]$. The $6$ elements of our set are the following ones.

\begin{footnotesize}
$$\begin{array}{|c|c|c|c|c|c|c|c|c|c|c|c|c|}\hline
\sha 30 &\sha  -12 &\sha  -18 &\shb  3 &\shb  -3 &\shc  5 &\shc  -5 &\shb  7 &\shb  -7 &\shc  9 &\shc  -9 &\shb  11 & \shb -11 \\ \hline
\sha -6 &\sha  85 &\sha  -79 &\shb  -161 &\shb  161 &\shc  -162 &\shc  162 &\shb  -163 &\shb  163 &\shc  -164 &\shc  164 & \shb  -165 &\shb  165 \\ \hline
\sha -24 &\sha  -73 &\sha  97 &\shb  158 &\shb  -158 &\shc  157 &\shc  -157 &\shb  156 &\shb  -156 &\shc  155 & \shc  -155 &\shb  154 &\shb  -154 \\ \hline
\shb 1 &\shb  -160 &\shb  159 &\shd  32 &\shd  -36 & \shd -2 & \shd 4 &\shd  -40 &\shd  42 &\sha  80 &\sha  -82 &\sha  -84 &\sha  86 \\ \hline
\shb -1 &\shb  160 &\shb  -159 &\shd  -32 &\shd  36 &\shd  2 &\shd  -4 &\shd  40 &\shd  -42 &\sha  -80 &\sha  82 &\sha  84 &\sha  -86 \\ \hline
 \end{array},$$
$$\begin{array}{|c|c|c|c|c|c|c|c|c|c|c|c|c|}\hline
 \sha 34 &\sha  -14 &\sha  -20 &\shb   15 &\shb   -15 &\shc 17 &\shc -17 & 19\shb  &\shb  -19 &\shc  21 &\shc  -21 &\shb  23 &\shb  -23 \\ \hline
 \sha -8 &\sha  89 &\sha  -81 &\shb   -167 &\shb   167 &\shc -168 &\shc 168 &\shb  -169 &\shb  169 &\shc  -170 &\shc  170 & \shb -171 &\shb  171 \\ \hline
\sha -26 &\sha  -75 &\sha  101 &\shb   152 & \shb  -152 & \shc 151 &\shc  -151 & \shb 150 &\shb  -150 & \shc 149 &\shc  -149 &\shb  148 &\shb  -148 \\ \hline
\shb  13 &\shb   -166 &\shb   153 & \shd  87 &\shd  -91 &\shd  -44 &\shd  46 &\shd  -48 &\shd  50 &\sha  88 &\sha  -90 &\sha  -92 &\sha  94 \\ \hline
\shb  -13 & \shb  166 &\shb   -153 &\shd  -87 &\shd  91 &\shd  44 &\shd  -46 & \shd 48 &\shd  -50 & \sha -88 & \sha 90 &\sha  92 &\sha  -94 \\ \hline
\end{array},$$
$$\begin{array}{|c|c|c|c|c|c|c|c|c|c|c|c|c|}\hline
 \sha 38 &\sha  -16 &\sha  -22 &\shb   27 &\shb   -27 &\shc  29 &\shc  -29 &\shb  31 &\shb  -31 &\shc  33 &\shc  -33 &\shb  35 &\shb  -35 \\ \hline
 \sha -10 &\sha  93 &\sha  -83 &\shb   -173 &\shb   173 & \shc -174 & \shc 174 &\shb  -175 &\shb  175 &\shc  -176 &\shc  176 &\shb  -177 & \shb 177 \\ \hline
 \sha -28 &\sha  -77 &\sha  105 &\shb   146 &\shb   -146 & \shc 145 & \shc -145 &\shb  144 &\shb  -144 & \shc 143 &\shc  -143 &\shb  142 & \shb -142 \\ \hline
\shb  25 &\shb   -172 &\shb   147 &\shd  95 & \shd -99 & \shd -52 &\shd  54 & \shd -56 &\shd  58 &\sha  96 &\sha  -98 & \sha -100 & \sha 102 \\ \hline
\shb  -25 &\shb   172 &\shb   -147 & \shd -95 &\shd  99 & \shd 52 &\shd  -54 & \shd 56 &\shd  -58 &\sha  -96  & \sha 98 &\sha  100 &\sha  -102 \\ \hline
\end{array},$$
 $$\begin{array}{|c|c|c|c|c|c|c|c|c|c|c|c|c|}\hline
 \sha -30 &\sha  12 &\sha  18 &\shb   39 &\shb   -39 & \shc 41 & \shc -41 &\shb  43 &\shb  -43 &\shc  45 &\shc  -45 &\shb  47 &\shb  -47 \\ \hline
\sha 6 &\sha  -85 &\sha  79 &\shb   -179 &\shb   179 & \shc -180 &\shc  180 &\shb  -181 &\shb  181 & \shc -182 & \shc 182 &\shb  -183 &\shb  183 \\ \hline
\sha 24 &\sha  73 &\sha  -97 &\shb   140 &\shb   -140 & \shc 139 &\shc  -139 & \shb 138 &\shb  -138 & \shc 137 &\shc  -137 &\shb  136 & \shb -136 \\ \hline
 \shb  37 &\shb   -178 &\shb   141 &\shd  106 &\shd  -110 &\shd  -60 &\shd  62 &\shd  -64 &\shd  66 & \sha 112 &\sha  -114 &\sha  -116 &\sha  118 \\ \hline
\shb  -37 &\shb   178 &\shb   -141 & \shd -106 &\shd  110 &\shd  60 & \shd -62 &\shd  64 &\shd  -66 & \sha -112 &\sha  114 &\sha  116 &\sha  -118 \\ \hline
\end{array},$$
 $$\begin{array}{|c|c|c|c|c|c|c|c|c|c|c|c|c|}\hline
 \sha -34 &\sha  14 &\sha  20 &\shb   51 & \shb  -51 &\shc  53 &\shc  -53 &\shb  55 &\shb  -55 &\shc  57 & \shc -57 &\shb  59 &\shb  -59 \\ \hline
\sha 8 &\sha  -89 &\sha  81 & \shb  -185 &\shb   185 &\shc  -186 & \shc 186 &\shb  -187 &\shb  187 &\shc  -188 & \shc 188 &\shb  -189 &\shb  189 \\ \hline
\sha 26 &\sha  75 &\sha  -101 &\shb   134 &\shb   -134 &\shc  133 &\shc  -133 &\shb  132 &\shb  -132 & \shc 131 & \shc -131 &\shb  130 &\shb  -130 \\ \hline
\shb  49 &\shb   -184 &\shb   135 & \shd 107 &\shd  -111 &\shd  -68 & \shd 70 &\shd  -72 &\shd  74 &\sha  120 &\sha  -122 &\sha  -113 & \sha 115 \\ \hline
\shb  -49 &\shb   184 &\shb   -135 &\shd  -107 &\shd  111 &\shd  68 &\shd  -70 & \shd 72 &\shd  -74 &\sha  -120 &\sha  122 &\sha  113 &\sha  -115 \\ \hline
 \end{array},$$
$$\begin{array}{|c|c|c|c|c|c|c|c|c|c|c|c|c|}\hline
\sha -38 &\sha  16 &\sha  22 & \shb  63 &\shb   -63 &\shc  65 &\shc  -65 &\shb  67 &\shb  -67 & \shc 69 &\shc  -69 &\shb  71 &\shb  -71 \\ \hline
\sha 10 &\sha  -93 &\sha  83 &\shb   -191 & \shb  191 &\shc  -192 &\shc  192 & \shb -193 &\shb  193 & \shc -194 &\shc  194 &\shb  -195 &\shb  195 \\ \hline
\sha 28 &\sha  77 &\sha  -105 &\shb   128 &\shb   -128 &\shc  127 &\shc  -127 &\shb  126 &\shb  -126 & \shc 125 &\shc  -125 &\shb  124 &\shb  -124 \\ \hline
\shb  61 &\shb   -190 &\shb   129 & \shd 76 &\shd  -78 & \shd -103 & \shd 104 & \shd -108 & \shd 109 &\sha  117 &\sha  -119 & \sha -121 &\sha  123 \\ \hline
 \shb  -61 &\shb   190 &\shb   -129 &\shd  -76 &\shd  78 &\shd  103 &\shd  -104 &\shd  108 &\shd  -109 &\sha  -117 &\sha  119 & \sha 121 &\sha  -123 \\ \hline
\end{array}.$$
 \end{footnotesize}
 \end{ex}

\begin{lem}\label{511-3}
Suppose that $a\geqslant 5$ and $b\geqslant 11$ are such that $a\equiv 1 \pmod 4$ and $b\equiv 3 \pmod 4$.
There exists an $\SMAS(a,b;2c)$ for all $c\geqslant 3$ such that $c\equiv 3 \pmod 4$.
\end{lem}

\begin{proof}
Write $c=4t+3$ where $t\geqslant 0$.
Our $\SMAS(a,b;2c)$ will consist
of $8t+4$ arrays of type 
\begin{equation}\label{97}
\begin{array}{|c|c|c|c|c|c|c|c|}\hline
  X_3      &    X_{2}^\T  & X_{2}^\T 
  &  X_2^\T & X_{2}^\T  & \cdots & X_2^\T & X_{2}^\T \\ \hline
  X_{2}    &  \multicolumn{2}{c|}{Z_4}   &  \multicolumn{2}{c|}{Z_4} & \cdots &  \multicolumn{2}{c|}{Z_4}    \\ \hline  
   X_{2}    &  \multicolumn{2}{c|}{Z_4}   &  \multicolumn{2}{c|}{Z_4} & \cdots &  \multicolumn{2}{c|}{Z_4}    \\ \hline  
  \vdots    &  \multicolumn{2}{c|}{\vdots}   &  \multicolumn{2}{c|}{\vdots} & \ddots &  \multicolumn{2}{c|}{\vdots}    \\ \hline  
   X_{2}    &  \multicolumn{2}{c|}{Z_4}   &  \multicolumn{2}{c|}{Z_4} & \cdots &  \multicolumn{2}{c|}{Z_4}    \\ \hline  
\end{array}
\end{equation}
and two arrays of type  
\begin{equation}\label{97special}
\begin{array}{|c|c|c|c|c|c|c|c|}\hline
\multicolumn{2}{|c|}{}  &    X_{2}^\T & X_{2}^\T  & X_{2}^\T & \cdots  & X_{2}^\T  & X_{2}^\T  
\\ \cline{3-8}
\multicolumn{2}{|c|}{\;\;\smash{\raisebox{.5\normalbaselineskip}{$Y_5$}}\;\;\;} & \multicolumn{1}{c|}{W_2} 
& \multicolumn{2}{c|}{Z_4}  & \cdots & \multicolumn{2}{c|}{Z_4}  \\ \hline  
   X_{2}    &  \multicolumn{2}{c|}{Z_4}   &  \multicolumn{2}{c|}{Z_4} & \cdots &  \multicolumn{2}{c|}{Z_4}    \\ \hline  
  \vdots    &  \multicolumn{2}{c|}{\vdots}   &  \multicolumn{2}{c|}{\vdots} & \ddots &  \multicolumn{2}{c|}{\vdots}    \\ \hline  
   X_{2}    &  \multicolumn{2}{c|}{Z_4}   &  \multicolumn{2}{c|}{Z_4} & \cdots &  \multicolumn{2}{c|}{Z_4}    \\ \hline  
\end{array},
\end{equation}
where each $X_r$ is a zero-sum block of size $r\times 3$, each $Z_4$ is a zero-sum block of size $2\times 4$,
each $Y_5$ is a block of size $5\times 5$ and each $W_2$ is a block of size $2\times 2$ such that
\begin{equation}\label{Y5}
\begin{array}{|c|c|c|}\hline
\multicolumn{2}{|c|}{}  &    X_{2}^\T  \\ \cline{3-3}
\multicolumn{2}{|c|}{\;\;\smash{\raisebox{.5\normalbaselineskip}{$Y_5$}}\;\;\;} & \multicolumn{1}{c|}{W_2}  \\ \hline  
\end{array} 
\end{equation}
is a zero-sum block of size $5\times 7$.

Write $a=4v+5$ and  $b=4w+11$, where $v,w\geqslant 0$.
Hence, we will construct
$8t+4$ blocks $X_3$, $4(v+w)(4t+3) +  40t+26$
blocks $X_2$, $(8t+4)(w+2)(2v+1)  + 2(w+1 + 2v(w+2))$ blocks  $Z_4$, 
$2$ blocks $Y_5$ and $2$ blocks $W_2$.
So, take
$$\begin{array}{rcl}
\mathfrak{A} & =& \mathfrak{A}(4t+5,\; 4(v+w)(4t+3)+ 40t+30,\; 4t+2 ),\\
\mathfrak{B} & =&  \mathfrak{B}(9 ,\; 8(2vw+ 5v+2w) (4t+3)+ 180t+131,\;  4(v+w)(4t+3) +  40t+26  ),\\
\mathfrak{C} & =& \mathfrak{C}(8(v+w)(4t+3)+140t+92  ,\; 8(2vw+5v+2w)(4t+3)+ 180t+132,\;\\
&&  8(v+w)(4t+3)+140t+102)
  \end{array}$$
 as in Lemmas \ref{BlockA},  \ref{BlockB} and
 \ref{BlockC}, respectively.
 We have $$\E(\mathfrak{A}) \dcup\E(\mathfrak{B})\dcup \E(\mathfrak{C})=\pm [1, 4(4vw+11v+5w)(4t+3) + 220t+165]
 \setminus \pm \left(\mathcal{M}_1\dcup \ldots \dcup \mathcal{M}_{12} \right),$$
 where
$$\begin{array}{rcl}
\mathcal{M}_1 & =&  [12, 8t+10]_2, \\
\mathcal{M}_2 & =& [40t+28, 40t+34]_2, \\
\mathcal{M}_3 & =& [40t+38,  56t+42]_4,  \\
\mathcal{M}_4 & =& [56t+44, 8(v+w) (4t+3)+  96t+66 ]_2, \\
\mathcal{M}_5 & =& [ 8(v+w)(4t+3)+96t+68,  8(v+w)(4t+3)+96t+75], \\
\mathcal{M}_6 & =& [8(v+w)(4t+3)+96t+76, 8(v+w)(4t+3)+128t+90]_2, \\
\mathcal{M}_7 & =& [ 8(v+w)(4t+3)+96t+79, 8(v+w)(4t+3)+128t+91 ]_4, \\
\mathcal{M}_8 & =& [8(v+w)(4t+3)+128t+92, 8(v+w)(4t+3)+140t+91], \\
\mathcal{M}_{9} & =& 
[8(v+w)(4t+3)+140t+95, 8(v+w)(4t+3)+140t+97]_2,\\   
\mathcal{M}_{10} & =& [8(v+w)(4t+3)+140t+106, 4(4vw+9v+3w)(4t+3)+140t+105  ],\\
\mathcal{M}_{11} & =& [8(v+w)(4t+3)+140t+93, 8(v+w)(4t+3)+140t+94],\\
\mathcal{M}_{12} & =&
[8(v+w)(4t+3)+140t+99, 8(v+w)(4t+3)+140t+100].
\end{array}$$
So, we replace the $8t+4$ instances of $X_3$ with the arrays in  $\mathfrak{A}$,  
the $2$  instances of $Y_5$ with the arrays in  $\mathfrak{C}$ and the $ 4(v+w)(4t+3) +  40t+26  $ instances of $X_2$ with the arrays in $\mathfrak{B}$.

Define 
$$\mathfrak{R}=\{U_2(\mathcal{M}_{11}), U_2(\mathcal{M}_{12}) \}.$$
Then, $\E(\mathfrak{R})=\pm (\mathcal{M}_{11}\dcup\mathcal{M}_{12})$: we replace the $2$ instances of $W_2$ with the arrays
in $\mathfrak{R}$. 
Next, write
$$\mathcal{M}_5=\mathcal{M}_5^1\dcup \mathcal{M}_5^2\dcup \mathcal{M}_5^3\dcup \mathcal{M}_5^4,$$
where
$$\begin{array}{rcl}
\mathcal{M}_5^1 & = & [8(v+w)(4t+3)+96t+68, 8(v+w)(4t+3)+96t+72]_4 , \\
\mathcal{M}_5^2 & = & [8(v+w)(4t+3)+96t+69, 8(v+w)(4t+3)+96t+70], \\
\mathcal{M}_5^3 & = & [8(v+w)(4t+3)+96t+71, 8(v+w)(4t+3)+96t+73]_2, \\
\mathcal{M}_5^4 & = & [8(v+w)(4t+3)+96t+74, 8(v+w)(4t+3)+96t+75].
 \end{array}$$
The set $\mathcal{M}_3 \dcup \mathcal{M}_5^1 \dcup\mathcal{M}_7$  can be written as a disjoint union 
$F_1\dcup \ldots \dcup F_{2\alpha}$, where $\alpha=3t+2$ and each $F_i$ is a $2$-set of type $4$.
The set $\mathcal{M}_1\dcup \mathcal{M}_2\dcup \mathcal{M}_4\dcup \mathcal{M}_5^3 \dcup 
 \mathcal{M}_6\dcup \mathcal{M}_{9}$  can be written as a disjoint union 
 $G_1\dcup \ldots \dcup G_{2\beta}$, where $\beta=(v+w)(4t+3)+10t+7 $  and each $G_j$ is a $2$-set of type $2$.
 Finally, the set $\mathcal{M}_5^2\dcup \mathcal{M}_5^4\dcup \mathcal{M}_8\dcup \mathcal{M}_{10}$
 can be written as a disjoint union  $H_1\dcup \ldots \dcup H_{2\gamma}$, where $\gamma= (4vw+7v+w) (4t+3)+  3t+1$ and each $H_h$ is a $2$-set of type $1$.
Call
 $$\begin{array}{rcl}
\mathfrak{S}_1 & =& \left\{Q_3(F_{2i-1},F_{2i} ): i \in [1,3t+2] \right\},\\
\mathfrak{S}_2 & =& \left\{Q_2(G_{2j-1},G_{2j} ): j \in [1,(v+w)(4t+3)+10t+7] \right\},\\
\mathfrak{S}_3 & =& \left\{Q_1(H_{2h-1},H_{2h} ): h \in [1,(4vw+7v+w)(4t+3)+3t+1] \right\}.
 \end{array}$$
 Then $\mathfrak{S}=\mathfrak{S}_1\dcup\mathfrak{S}_2\dcup \mathfrak{S}_3 $ consists of
 $2(2vw+4v+w)(4t+3)+ 16t+10 $ zero-sum blocks of size $2\times 4$ and
 $\E(\mathfrak{S})=\pm  (\mathcal{M}_1\dcup \ldots\dcup \mathcal{M}_{10})$.
Hence, we replace the instances of $Z_4$ with the arrays in~$\mathfrak{S}$.
\end{proof}

\begin{lem}\label{97-3}
Suppose that $a\geqslant 9$ is such that  $a\equiv 1 \pmod 4$.
There exists an $\SMAS(a,7;2c)$ for all $c\geqslant 3$ such that $c\equiv 3 \pmod 4$.
\end{lem}

\begin{proof}
Write $c=4t+3$ where $t\geqslant 0$.
Our $\SMAS(a,7;2c)$ will consist
of $8t+4$ arrays of type 
\begin{equation}\label{97sm}
\begin{array}{|c|c|c|}\hline
  X_3      &    X_{2}^\T  & X_{2}^\T \\ \hline
  X_{2}    &  \multicolumn{2}{c|}{Z_4}    \\ \hline  
   X_{2}    &  \multicolumn{2}{c|}{Z_4}     \\ \hline  
  \vdots    &  \multicolumn{2}{c|}{\vdots}   \\ \hline  
   X_{2}    &  \multicolumn{2}{c|}{Z_4}      \\ \hline  
\end{array}
\end{equation}
and two arrays of type  
\begin{equation}\label{97supsp}
\begin{array}{|c|c|c|}\hline
\multicolumn{2}{|c|}{}  &    X_{2}^\T \\ \cline{3-3}
\multicolumn{2}{|c|}{\;\;\smash{\raisebox{.5\normalbaselineskip}{$Y_5$}}\;\;\;} & \multicolumn{1}{c|}{W_2}   \\ \hline  
   X_{2}    &  \multicolumn{2}{c|}{Z_4}     \\ \hline  
  \vdots    &  \multicolumn{2}{c|}{\vdots}   \\ \hline  
   X_{2}    &  \multicolumn{2}{c|}{Z_4}    \\ \hline  
\end{array},
\end{equation}
where each $X_r$ is a zero-sum block of size $r\times 3$, each $Z_4$ is a zero-sum block of size $2\times 4$,
each $Y_5$ is a block of size $5\times 5$ and each $W_2$ is a block of size $2\times 2$ such that \eqref{Y5} is a zero-sum block of size $5\times 7$.

Write $a=4v+9$, where $v\geqslant 0$.
Hence, we will construct
$8t+4$ blocks $X_3$, $4v(4t+3)+40t+26$
blocks $X_2$, $(8t+4) (2v+3)+ 2(2v+2)$ blocks  $Z_4$, 
$2$ blocks $Y_5$ and $2$ blocks $W_2$.
So, take
$$\begin{array}{rcl}
\mathfrak{A} & =& \mathfrak{A}(4t+5,\; 4v(4t+3)+40t+30 ,\; 4t+2 ),\\
\mathfrak{B} & =&  \mathfrak{B}(9 ,\;  24v(4t+3)+212t+155  ,\;  4v(4t+3)+40t+26 ),\\
\mathfrak{C} & =& \mathfrak{C}( 8v(4t+3)+128t+92,\;    24v(4t+3)+212t+156,\;8v(4t+3)+128t+102)
  \end{array}$$
 as in Lemmas \ref{BlockA},  \ref{BlockB} and \ref{BlockC}, respectively.
 We have $$\E(\mathfrak{A}) \dcup\E(\mathfrak{B})\dcup \E(\mathfrak{C})=\pm [1, 28v(4t+3)+ 252t+189]
 \setminus \pm \left(\mathcal{M}_1\dcup \ldots \dcup \mathcal{M}_{12} \right),$$
 where 
$$\begin{array}{rcl}
\mathcal{M}_1 & =& [12, 8t+10]_2, \\
\mathcal{M}_2 & =&  [40t+28, 40t+34]_2, \\
\mathcal{M}_3 & =&  [40t+38, 56t+42]_4,\\
\mathcal{M}_4 & =& [56t+44, 8v(4t+3)+96t+66]_2, \\
\mathcal{M}_5 & =&  [8v(4t+3)+ 96t+68, 8v(4t+3)+96t+75], \\
\mathcal{M}_6 & =&  [ 8v(4t+3)+96t+76,8v(4t+3)+128t+90]_2, \\
\mathcal{M}_7 & =& [8v(4t+3)+96t+79, 8v(4t+3)+128t+91]_4,  \\
\mathcal{M}_8 & =& [ 8v(4t+3) +128t+106,  8v(4t+3)+172t+129],\\
\mathcal{M}_{9} & =& [ 8v(4t+3) +128t+93,  8v(4t+3)+128t+94],
\\   
\mathcal{M}_{10} & =&[ 8v(4t+3)+128t+95,8v(4t+3)+128t+97]_2,
\\
\mathcal{M}_{11} & =& [ 8v(4t+3)+128t+99,8v(4t+3)+128t+100],\\
\mathcal{M}_{12} & =&[8v(4t+3)+172t+130,  20v(4t+3)+172t+129].
\end{array}$$
So, we replace the $8t+4$ instances of $X_3$ with the arrays in  $\mathfrak{A}$,  
the $2$  instances of $Y_5$ with the arrays in  $\mathfrak{C}$ and the $4v(4t+3)+40t+26$ instances of $X_2$ with the arrays in $\mathfrak{B}$.

Define 
$$\mathfrak{R}=\{U_2(\mathcal{M}_{9}), U_2(\mathcal{M}_{11}) \}.$$
Then, $\E(\mathfrak{R})=\pm(\mathcal{M}_{9}\dcup\mathcal{M}_{11})$: we replace the $2$ instances of $W_2$ with the arrays
in $\mathfrak{R}$. 
Next, write
$$\mathcal{M}_5=\mathcal{M}_5^1\dcup \mathcal{M}_5^2\dcup \mathcal{M}_5^3\dcup \mathcal{M}_5^4,$$
where
$$\begin{array}{rcl}
\mathcal{M}_5^1 & = & [8v(4t+3)+96t+68, 8v(4t+3)+96t+72]_4 , \\
\mathcal{M}_5^2 & = & [8v(4t+3)+96t+69, 8v(4t+3)+96t+70], \\
\mathcal{M}_5^3 & = & [8v(4t+3)+96t+71, 8v(4t+3)+96t+73]_2, \\
\mathcal{M}_5^4 & = & [8v(4t+3)+96t+74, 8v(4t+3)+96t+75].
 \end{array}$$
The set $\mathcal{M}_3 \dcup \mathcal{M}_5^1 \dcup\mathcal{M}_7$  can be written as a disjoint union 
$F_1\dcup \ldots \dcup F_{2\alpha}$, where $\alpha=3t+2$ and each $F_i$ is a $2$-set of type $4$.
The set $\mathcal{M}_1\dcup \mathcal{M}_2\dcup \mathcal{M}_4\dcup \mathcal{M}_5^3 \dcup 
 \mathcal{M}_6\dcup \mathcal{M}_{10}$  can be written as a disjoint union 
 $G_1\dcup \ldots \dcup G_{2\beta}$, where $\beta=v(4t+3)+10t+7 $  and each $G_j$ is a $2$-set of type $2$.
 Finally, the set $\mathcal{M}_5^2\dcup \mathcal{M}_5^4\dcup \mathcal{M}_8\dcup \mathcal{M}_{12}$ 
 can be written as a disjoint union  $H_1\dcup \ldots \dcup H_{2\gamma}$, where $\gamma=  3v(4t+3)+11t+7$ and each $H_h$ is a $2$-set of type $1$.
Call
 $$\begin{array}{rcl}
\mathfrak{S}_1 & =& \left\{Q_3(F_{2i-1},F_{2i} ): i \in [1,3t+2] \right\},\\
\mathfrak{S}_2 & =& \left\{Q_2(G_{2j-1},G_{2j} ): j \in [1,v(4t+3)+10t+7] \right\},\\
\mathfrak{S}_2 & =& \left\{Q_1(H_{2h-1},H_{2h} ): h \in [1,3v(4t+3)+11t+7] \right\}.
 \end{array}$$
 Then $\mathfrak{S}=\mathfrak{S}_1\dcup\mathfrak{S}_2\dcup \mathfrak{S}_3 $ consists of
 $ 4v(4t+3)+ 24t+16 $ zero-sum blocks of size $2\times 4$ and
 $\E(\mathfrak{S})=\pm  (\mathcal{M}_1\dcup \ldots\dcup \mathcal{M}_{8}\dcup \mathcal{M}_{10}\dcup \mathcal{M}_{12})$.
Hence, we replace the instances of $Z_4$ with the arrays in  $\mathfrak{S}$.
\end{proof}

\subsection{Subcase \texorpdfstring{$c\equiv 1 \pmod 4$}{c ≡ 1 (mod 4)}}

\begin{lem}\label{59-1}
Suppose that $a\geqslant 5$ and $b\geqslant 9$ are such that  $a\equiv b\equiv 1 \pmod 4$.
There exists an $\SMAS(a,b;2)$.
\end{lem}

\begin{proof}
Our $\SMAS(a,b;2)$ will consist of two arrays of type  
\begin{equation}\label{extra9}
\begin{array}{|c|c|c|c|c|c|c|c|c|}\hline
\multicolumn{2}{|c|}{}  &    X_{2}^\T &    X_{2}^\T & X_{2}^\T  & X_{2}^\T & \cdots  & X_{2}^\T  & X_{2}^\T  
\\ \cline{3-9}
\multicolumn{2}{|c|}{\;\;\smash{\raisebox{.5\normalbaselineskip}{$Y_5$}}\;\;\;} & \multicolumn{2}{c|}{W_4} 
& \multicolumn{2}{c|}{Z_4}  & \cdots & \multicolumn{2}{c|}{Z_4}  \\ \hline  
   X_{2}    &  & \multicolumn{2}{c|}{Z_4}   &  \multicolumn{2}{c|}{Z_4} & \cdots &  \multicolumn{2}{c|}{Z_4}    \\   \cline{1-1} \cline{3-9}
   X_{2}    & \smash{\raisebox{.5\normalbaselineskip}{$Z_4^\T$}}  &\multicolumn{2}{c|}{Z_4}   &  \multicolumn{2}{c|}{Z_4} & \cdots &  \multicolumn{2}{c|}{Z_4}    \\ \hline  
   \vdots    & \vdots &  \multicolumn{2}{c|}{\vdots}   &  \multicolumn{2}{c|}{\vdots} & \ddots &  \multicolumn{2}{c|}{\vdots}    \\ \hline  
      X_{2}    &  & \multicolumn{2}{c|}{Z_4}   &  \multicolumn{2}{c|}{Z_4} & \cdots &  \multicolumn{2}{c|}{Z_4}    \\   \cline{1-1} \cline{3-9}
   X_{2}    & \smash{\raisebox{.5\normalbaselineskip}{$Z_4^\T$}}  &\multicolumn{2}{c|}{Z_4}   &  \multicolumn{2}{c|}{Z_4} & \cdots &  \multicolumn{2}{c|}{Z_4}    \\ \hline  
\end{array},
\end{equation}
where  each $X_2$ is a zero-sum block of size $2\times 3$,
each $Z_4$ is a zero-sum block of size $2\times 4$,
each $Y_5$ is a block of size $5\times 5$
and each $W_4$ is a block of size $2\times 4$ such that
\begin{equation}\label{Y5b}
\begin{array}{|c|c|c|c|}\hline
 \multicolumn{2}{|c|}{}  &  X_{2}^\T  & X_{2}^\T  \\ \cline{3-4}
\multicolumn{2}{|c|}{\;\;\smash{\raisebox{.5\normalbaselineskip}{$Y_5$}}\;\;\;} & \multicolumn{2}{c|}{W_4}    \\ \hline  
\end{array} 
\end{equation}
is a zero-sum block of size $5\times 9$. 

Write $a=4v+5$ and $b=4w+9$, where $v,w\geqslant 0$.
Hence, we will construct
$4(v+w+1)$
blocks $X_2$, $2(3v+ (2v+1)w)$ blocks  $Z_4$, 
$2$ blocks $Y_5$ and $2$ blocks $W_4$.
So, take
$$\begin{array}{rcl}
 \mathfrak{B} & =&  \mathfrak{B}(9 ,\; 16vw+32v+16w+33 ,\; 4(v+w+1)),\\
 \mathfrak{C} & =& \mathfrak{C}(8(v+w)+16,\; 16vw+32v+16w+34 ,\; 8(v+w)+26)
 \end{array}$$
 as in Lemmas \ref{BlockB} and \ref{BlockC}, respectively.
 We have $$\E(\mathfrak{B})\dcup \E(\mathfrak{C})=\pm [1, 16vw+36v+20w+45 ]\setminus \pm \left(\mathcal{M}_1\dcup \ldots \dcup \mathcal{M}_{5} \right),$$
where
$$\begin{array}{rcl}
\mathcal{M}_1 & =& [12, 8(v+w)+14]_2,\\
\mathcal{M}_2 & =& [8(v+w)+17, 8(v+w)+18], \\
\mathcal{M}_3 & =& [8(v+w)+19, 8(v+w)+21]_2, \\
\mathcal{M}_4 & =& [8(v+w)+23, 8(v+w)+24],\\
\mathcal{M}_5 & =& [8(v+w)+30, 16v(w+1) + 12(v+w) +29].
\end{array}$$
So, we replace the $2$  instances of $Y_5$ with the arrays in  $\mathfrak{C}$ and the $ 4(v+w+1)$ instances of $X_2$ with the arrays in $\mathfrak{B}$.

Write 
$$\mathcal{M}_1 = \mathcal{M}_1^1 \dcup \mathcal{M}_1^2, \quad \text{ where } \quad  \mathcal{M}_1^1 =[12, 14]_2 \equad \mathcal{M}_1^2=[16, 8(v+w)+14]_2.$$
Define 
$$\mathfrak{R}=\{U_1(\mathcal{M}_1^1,\mathcal{M}_2), U_1(\mathcal{M}_3,\mathcal{M}_4) \}.$$
Then, $\E(\mathfrak{R})=\pm(\mathcal{M}_1^1\dcup \mathcal{M}_2\dcup \mathcal{M}_3\dcup \mathcal{M}_4)$: 
we replace the $2$ instances of $W_4$ with the arrays in~$\mathfrak{R}$. 

The set $\mathcal{M}_1^2$ can be written as a disjoint union $G_1\dcup \ldots \dcup G_{2\alpha}$, where $\alpha=v+w$ and each $G_{j}$ is a $2$-set of type $2$.
The set $\mathcal{M}_5$ can be written as a disjoint union $H_1\dcup \ldots \dcup H_{2\gamma}$, where $\gamma=4v(w+1)+(v+w)$ and each $H_{h}$ is a $2$-set of type $1$ .
Call
$$\begin{array}{rcl}
\mathfrak{S}_1 & =& \left\{Q_2(G_{2j-1},G_{2j} ): j \in [1,v+w] \right\},\\
\mathfrak{S}_2 & =& \left\{ Q_1(H_{2h-1},H_{2h}) : h \in [1,4v(w+1)+(v+w)] \right\}.
\end{array}$$
Then $\mathfrak{S}=\mathfrak{S}_1\dcup \mathfrak{S}_2$ consists of $4vw +6v+ 2w$ zero-sum blocks of size $2\times 4$ and
$\E(\mathfrak{S})=\pm (\mathcal{M}_1^2\dcup \mathcal{M}_{5})$. 
Hence, we replace the instances of $Z_4$ with the arrays in~$\mathfrak{S}$.
\end{proof}

\begin{lem}\label{a1b1-1}
Suppose that $a\geqslant 5$ and $b\geqslant 9$ are such that $a\equiv b\equiv 1 \pmod 4$.
There exists an $\SMAS(a,b;2c)$ for all $c\geqslant 1$ such that $c\equiv 1 \pmod 4$.
\end{lem}

\begin{proof}
By Lemma \ref{59-1} we can write $c=4t+5$ where $t\geqslant 0$. 
Our $\SMAS(a,b;2c)$ will consist
of $8t+8$ arrays of type
\begin{equation}\label{97b}
\begin{array}{|c|c|c|c|c|c|c|c|c|}\hline
  X_3      &   X_2^\T & X_{2}^\T  & X_{2}^\T 
  &  X_2^\T & X_{2}^\T  & \cdots & X_2^\T & X_{2}^\T \\ \hline
  X_{2}    &  \multicolumn{3}{c|}{Z_6} 
  &  \multicolumn{2}{c|}{Z_4} & \cdots &  \multicolumn{2}{c|}{Z_4} 
   \\ \hline  
      X_{2}    &  & \multicolumn{2}{c|}{Z_4}   &  \multicolumn{2}{c|}{Z_4} & \cdots &  \multicolumn{2}{c|}{Z_4}    \\   \cline{1-1} \cline{3-9}
   X_{2}    & \smash{\raisebox{.5\normalbaselineskip}{$Z_4^\T$}}  &\multicolumn{2}{c|}{Z_4}   &  \multicolumn{2}{c|}{Z_4} & \cdots &  \multicolumn{2}{c|}{Z_4}    \\ \hline  
   \vdots    & \vdots &  \multicolumn{2}{c|}{\vdots}   &  \multicolumn{2}{c|}{\vdots} & \ddots &  \multicolumn{2}{c|}{\vdots}    \\ \hline  
      X_{2}    &  & \multicolumn{2}{c|}{Z_4}   &  \multicolumn{2}{c|}{Z_4} & \cdots &  \multicolumn{2}{c|}{Z_4}    \\   \cline{1-1} \cline{3-9}
   X_{2}    & \smash{\raisebox{.5\normalbaselineskip}{$Z_4^\T$}}  &\multicolumn{2}{c|}{Z_4}   &  \multicolumn{2}{c|}{Z_4} & \cdots &  \multicolumn{2}{c|}{Z_4}    \\ \hline  
\end{array}
\end{equation}
and two arrays of type   \eqref{extra9},
where each $X_r$ is a zero-sum block of size $r\times 3$, each $Z_\ell$ is a zero-sum block of size $2\times \ell$,
each $Y_5$ is a block of size $5\times 5$ and each $W_4$ is a block of size $2\times 4$ such that \eqref{Y5b} is a zero-sum block of size $5\times 9$.

Write $a=4v+5$ and $b=4w+9$, where $v,w\geqslant 0$.
Hence, we will construct
$8t+8$ blocks $X_3$, $4(v+w)(4t+5)+    32t+36$
blocks $X_2$, $(8t+10) (3v+ w(2v+1))$ blocks  $Z_4$,  $8t+8$ blocks $Z_6$,
$2$ blocks $Y_5$ and $2$ blocks $W_4$.
So, take
$$\begin{array}{rcl}
\mathfrak{A} & =& \mathfrak{A}(4t+5,\; 4(v+w)(4t+5)+  32t+40   ,\;4t+4 ),\\
 \mathfrak{B} & =&  \mathfrak{B}(9 ,\; (16v(w+1)+16(v+w))(4t+5)+   148t+181 ,\;  4(v+w)(4t+5)+    32t+36  ),\\
 \mathfrak{C} & =& \mathfrak{C}(8(v+w)(4t+5)+  112t+130 ,\; ( 16v(w+1)+16(v+w))(4t+5)+ 148t+182 ,\; \\
 && 8(v+w)(4t+5)+ 112t+140 )
 \end{array}$$
 as in Lemmas \ref{BlockA},  \ref{BlockB} and \ref{BlockC}, respectively.
We have $$\E(\mathfrak{A}) \dcup\E(\mathfrak{B})\dcup\E( \mathfrak{C})=\pm [1, (16v(w+1)+20(v+w))(4t+5)+ 180t+225 ]
 \setminus \pm \left(\mathcal{M}_1\dcup \ldots \dcup \mathcal{M}_{12} \right),$$
where
$$\begin{array}{rcl}
\mathcal{M}_1 & =&  [12,8t+10]_2,  \\
\mathcal{M}_2 & =&   [40t+44, 40t+46]_2, \\
\mathcal{M}_3 & =&   [40t+50, 56t+62]_4, \\
\mathcal{M}_4 & =&  [56t+64,  8(v+w)(4t+5)+ 80t+94]_2 ,\\
\mathcal{M}_5 & =&  [8(v+w)(4t+5)+  80t+96, 8(v+w)(4t+5)+ 80t+97], \\
\mathcal{M}_6 & =&   [ 8(v+w)(4t+5)+ 80t+98,  8(v+w)(4t+5)+ 112t+128]_2 ,\\
\mathcal{M}_7 & =&  [ 8(v+w)(4t+5)+ 80t+99,   8(v+w)(4t+5)+ 112t+127 ]_4, \\
\mathcal{M}_8 & =&  [ 8(v+w)(4t+5)+ 112t+131, 8(v+w)(4t+5)+ 112t+135]_4,\\
\mathcal{M}_9 & =&  [ 8(v+w)(4t+5)+ 112t+132, 8(v+w)(4t+5)+ 112t+133], \\
\mathcal{M}_{10} & =& [ 8(v+w)(4t+5)+ 112t+137, 8(v+w)(4t+5)+ 112t+138],\\
\mathcal{M}_{11} & =&   [ 8(v+w)(4t+5)+ 112t+144, 8(v+w)(4t+5)+ 116t+145],\\
\mathcal{M}_{12} & =&  [8(v+w)(4t+5)+ 116t+146,  (16v(w+1)+12(v+w))(4t+5)+  116t+145].
\end{array}$$
So, we replace the $8t+8$ instances of $X_3$ with the arrays in  $\mathfrak{A}$,  
the $2$  instances of $Y_5$ with the arrays in  $\mathfrak{C}$ and the $4(v+w)(4t+5)+ 32t+36$ instances of $X_2$ with the arrays in $\mathfrak{B}$.

Write
$$\mathcal{M}_1=\mathcal{M}_1^1\dcup \mathcal{M}_1^2, \quad
\mathcal{M}_{6}=\mathcal{M}_{6}^1\dcup \mathcal{M}_{6}^2\dcup \mathcal{M}_{6}^3 \equad 
\mathcal{M}_{11}=\mathcal{M}_{11}^1\dcup \mathcal{M}_{11}^2\dcup \mathcal{M}_{11}^3,$$
where
$$\begin{array}{rcl}
\mathcal{M}_1^1 & =& [12, 8t+8]_4,\\
\mathcal{M}_1^2 & =& [14, 8t+10]_4,\\
\mathcal{M}_{6}^1 & =& [8(v+w)(4t+5)+80t+98, 8(v+w)(4t+5)+80t+100 ]_2,\\
\mathcal{M}_{6}^2 & =& [8(v+w)(4t+5)+80t+102, 8(v+w)(4t+5)+ 80t+104]_2,\\
\mathcal{M}_{6}^3 & =& [8(v+w)(4t+5)+80t+106, 8(v+w)(4t+5)+112t+128 ]_2,\\
\mathcal{M}_{11}^1 & =& [8(v+w)(4t+5)+ 112t+144, 8(v+w)(4t+5)+ 116t+142]_2,\\
\mathcal{M}_{11}^2 & =& [8(v+w)(4t+5)+ 112t+145, 8(v+w)(4t+5)+ 116t+143]_2,\\
\mathcal{M}_{11}^3 & =& [8(v+w)(4t+5)+ 116t+144, 8(v+w)(4t+5)+ 116t+145].
\end{array}$$

Define 
$$\mathfrak{R}=\{U_1( \mathcal{M}_6^1,
\mathcal{M}_5 ), U_1(\mathcal{M}_6^2 ,
\mathcal{M}_9 ) \}.$$
Then, $\E(\mathfrak{R})=\pm(\mathcal{M}_5\dcup \mathcal{M}_6^1\dcup \mathcal{M}_6^2\dcup \mathcal{M}_9)$: we replace the $2$ instances of $W_4$ with the arrays in~$\mathfrak{R}$. 

The set $\mathcal{M}_1^1\dcup \mathcal{M}_1^2\dcup \mathcal{M}_3 \dcup \mathcal{M}_7\dcup\mathcal{M}_8$  can be written as a disjoint union 
$F_1\dcup \ldots \dcup F_{\alpha}$, where $\alpha=8t+7$ and each $F_i$ is a $2$-set of type $4$.
The set $\mathcal{M}_2\dcup \mathcal{M}_4\dcup \mathcal{M}_{6}^3\dcup \mathcal{M}_{11}^1\dcup \mathcal{M}_{11}^2$ can be written as a disjoint union 
$G_1\dcup \ldots \dcup G_{\beta}$,  where $\beta=2(v+w)(4t+5)+16t+15$ and each $G_j$ is a $2$-set of type $2$.
Finally, the set  $\mathcal{M}_{10}\dcup \mathcal{M}_{11}^3\dcup \mathcal{M}_{12}$ can be written as a disjoint union 
$H_1\dcup \ldots \dcup H_{\gamma}$,  where $\gamma=(8v(w+1)+2(v+w))(4t+5)+2$ and each $H_h$ is a $2$-set of type $1$.
Call
$$\begin{array}{rcl}
\mathfrak{S}_1 & =& \left\{R_2(F_i,G_{2i-1},G_{2i} ): i \in [1,8t+7] \right\},\\
\mathfrak{S}_2 & =& \left\{R_1(G_{ \beta}, H_{\gamma-1} ,H_{\gamma} )  \right\},\\
\mathfrak{T}_1 & =& \left\{Q_2(G_{2j-1},G_{2j} ): j \in [8t+8, (v+w)(4t+5)+8t+7 ] \right\},\\
\mathfrak{T}_2 & =& \left\{Q_1(H_{2h-1},H_{2h} ): h \in [1, (4v(w+1)+(v+w))(4t+5) ] \right\}.
\end{array}$$
Then $\mathfrak{S}=\mathfrak{S}_1\dcup \mathfrak{S}_2$ consists of $8t+8$ zero-sum blocks of size $2\times 6$ and
the set $\mathfrak{T}=\mathfrak{T}_1\dcup \mathfrak{T}_2$ consists of $ (4v(w+1)+2(v+w))(4t+5)$ zero-sum blocks of size $2\times 4$.
Furthermore, $\E(\mathfrak{S})\dcup \E(\mathfrak{T})=\pm (\mathcal{M}_1\dcup \ldots \dcup \mathcal{M}_4 \dcup 
 \mathcal{M}_{6}^3\dcup \mathcal{M}_7\dcup \mathcal{M}_8\dcup  \mathcal{M}_{10} \dcup \mathcal{M}_{11} \dcup \mathcal{M}_{12})$.
Hence, we replace the instances of $Z_6$ with the arrays in  $\mathfrak{S}$
and the instances of $Z_4$ with the arrays in  $\mathfrak{T}$.
\end{proof}

\begin{lem}\label{511-1}
Suppose that $b\geqslant 11$ is such that $b\equiv 3 \pmod 4$. 
There exists an $\SMAS(5,b;2c)$ for all $c\geqslant 1$ such that $c\equiv 1 \pmod 4$.
\end{lem}

\begin{proof}
Write $c=4t+1$ where $t\geqslant 0$.
Our $\SMAS(5,b; 2c)$ will consist of $8t+2$ arrays of type 
\begin{equation}\label{57}
\begin{array}{|c|c|c|c|c|c|c|c|}\hline
  X_3      &    X_{2}^\T  & X_{2}^\T 
  &  X_2^\T & X_{2}^\T  & \cdots & X_2^\T & X_{2}^\T \\ \hline
  X_{2}    &  \multicolumn{2}{c|}{Z_4} 
  &  \multicolumn{2}{c|}{Z_4} & \cdots &  \multicolumn{2}{c|}{Z_4} 
   \\ \hline  
\end{array}.
\end{equation}

Write $b=4w+11$, where $w\geqslant 0$.
Hence, we will construct
$8t+2$ blocks $X_3$, $2(2w+5)(4t+1)$
blocks $X_2$ and $2(w+2)(4t+1)$ blocks  $Z_4$.
So, take
$$
\begin{array}{rcl}
\mathfrak{A} & =& \mathfrak{A}(4t,\;  2(2w+5)(4t+1)    ,\; 4t+1),\\
\mathfrak{B} & =& \mathfrak{B}(1,\;   (16w+45)(4t+1)   , \;  2(2w+5)(4t+1) )
\end{array}$$
as in Lemmas \ref{BlockA} and \ref{BlockB}, respectively. Then, 
$$\E(\mathfrak{A})\dcup \E(\mathfrak{B})=\pm [1,5(4w+11)(4t+1)]\setminus \pm \left(\mathcal{M}_1\dcup \ldots \dcup \mathcal{M}_5 \right),$$
where
$$\begin{array}{rcl}
\mathcal{M}_1 & =& [2,8t]_2, \\
\mathcal{M}_2 & =& [40t+12, 56t+8]_4, \\
\mathcal{M}_3 & =& [56t+12,  8w(4t+1)+ 128t+30]_2,\\
\mathcal{M}_4 & =& [8w(4t+1)+  96t+27,  8w(4t+1)+ 128t+31]_4, \\
\mathcal{M}_5 & =& [8w(4t+1)+ 128t+32, 12w(4t+1)+ 140t+35].
\end{array}$$
So, we replace the $8t+2$ instances of $X_3$ with the arrays in  $\mathfrak{A}$ and the $2(2w+5)(4t+1)$ instances of $X_2$ with the arrays in $\mathfrak{B}$.

Next, write
$$\mathcal{M}_3=\mathcal{M}_3^1\dcup \mathcal{M}_3^2
\equad 
\mathcal{M}_4=\mathcal{M}_4^1\dcup \mathcal{M}_4^2,
$$
where
$$\begin{array}{rcl}
\mathcal{M}_3^1 & =& [56t+12,  8w(4t+1)+128t+26]_2, \\
\mathcal{M}_3^2 & =& [8w(4t+1)+ 128t+28,  8w(4t+1)+ 128t+30]_2, \\
\mathcal{M}_4^1 & =& [8w(4t+1)+ 96t+27, 8w(4t+1)+128t+23 ]_4, \\
\mathcal{M}_4^2 & =& [8w(4t+1)+ 128t+27,  8w(4t+1)+ 128t+31]_4.
\end{array}$$
Note that 
$\mathcal{M}_3^2 \dcup
\mathcal{M}_4^2 = 
[8w(4t+1)+ 128t+27, 8w(4t+1)+ 128t+28]
\dcup [8w(4t+1)+ 128t+30,  8w(4t+1)+ 128t+31]$.

The set $\mathcal{M}_2\dcup \mathcal{M}_4^1$  can be written as a disjoint union 
$F_1\dcup \ldots \dcup F_{6t}$, where each $F_i$ is
a $2$-set of type $4$.
The set $\mathcal{M}_1\dcup \mathcal{M}_3^1$
can be written as a disjoint union 
$G_1\dcup \ldots \dcup G_{2\beta}$, where $\beta= (w+2)(4t+1)+2t $ and each $G_j$ is
a $2$-set of type $2$.
Finally, the set $\mathcal{M}_3^2\dcup \mathcal{M}_4^2\dcup \mathcal{M}_5$
can be written as a disjoint union 
$H_1\dcup \ldots \dcup H_{2\gamma}$, where $\gamma=w(4t+1)+ 3t+2 $ and each $H_h$ is
a $2$-set of type $1$.
Call
$$\begin{array}{rcl}
\mathfrak{R}_1 & =& \left\{Q_3(F_{2i-1},F_{2i} ):
 i \in [1,3t] \right\},\\
\mathfrak{R}_2 & =& \left\{Q_2(G_{2j-1},G_{2j} ):
 j \in [1, (w+2)(4t+1)+2t] \right\},\\
\mathfrak{R}_3 & =& \left\{Q_1(H_{2h-1},H_{2h} ):
 h \in [1,w(4t+1)+ 3t+2] \right\}.
\end{array}$$
Then $\mathfrak{R}=\mathfrak{R}_1\dcup \mathfrak{R}_2\dcup 
\mathfrak{R}_3$ consists of $2(w+2)(4t+1)$ zero-sum blocks of size $2\times 4$ and
$\E(\mathfrak{R})=\pm  (\mathcal{M}_1\dcup \ldots\dcup \mathcal{M}_5)$.
Hence, we replace the instances of $Z_4$ with the arrays in  $\mathfrak{R}$.
\end{proof}

\begin{lem}\label{77-1}
Suppose that $a,b\geqslant 7$ are such that $a\equiv b\equiv 3 \pmod 4$.
There exists an $\SMAS(a,b;2)$.
\end{lem}

\begin{proof}
Our $\SMAS(a,b;2)$ will consist of
two arrays of type \eqref{97special}.
Write $a=4v+7$ and $b=4w+7$, where $v,w\geqslant 0$.
Hence, we will construct
$4(v+w+1)$ blocks $X_2$,
$2(w+(2v+1)(w+1))$ blocks  $Z_4$,
$2$ blocks $Y_5$ and $2$ blocks $W_2$.
So, take
$$\begin{array}{rcl}
  \mathfrak{B} & =& \mathfrak{B}(9 ,\; 
  16vw + 24(v+w)+37, \; 4(v+w+1)    ) ,\\
  \mathfrak{C} & =& \mathfrak{C}(8(v+w) + 16  ,\;  16vw+24(v+w)+38, \; 8(v+w)+26)
\end{array}$$
as in Lemmas \ref{BlockB} and  \ref{BlockC}, respectively.
We have $$\E(\mathfrak{B})\dcup \E(\mathfrak{C})=\pm [1, 16vw+ 28(v+w)+49]\setminus \pm \left(\mathcal{M}_1\dcup \ldots \dcup \mathcal{M}_{5} \right),$$
where
$$\begin{array}{rcl}
\mathcal{M}_1 & =& [12, 8(v+w)+14]_2,\\
\mathcal{M}_2 & =& [8(v+w)+17, 8(v+w)+18],\\
\mathcal{M}_3 & =& [8(v+w)+19, 8(v+w)+21]_2,\\
\mathcal{M}_4 & =& [8(v+w)+23, 8(v+w)+24],\\
\mathcal{M}_5 & =& [8(v+w)+30, 
16vw+20(v+w)+ 33].
\end{array}$$
So, we replace the $2$ instances of $Y_5$ with the arrays in  $\mathfrak{C}$ and the $4(v+w+1)$ 
instances of $X_2$ with the arrays in $\mathfrak{B}$.
Call
$$\mathfrak{R} =\{U_2(\mathcal{M}_2), U_2(\mathcal{M}_4) \}.$$
Then, $\E(\mathfrak{R})=\pm \left(\mathcal{M}_2\dcup \mathcal{M}_4\right)$: we replace the $2$ instances of $W_2$ with the arrays in $\mathfrak{R}$.

The set $\mathcal{M}_1\dcup \mathcal{M}_3$ can be written as a disjoint union 
$G_1\dcup \ldots \dcup G_{2\beta}$, where $\beta=v+w+1$ and each $G_j$ is a $2$-set of type $2$.
The set $\mathcal{M}_{5}$ can be written as a disjoint union 
$H_1\dcup \ldots \dcup H_{2\gamma}$, where $\gamma=4vw+3(v+w)+1$ and each $H_h$ is a $2$-set of type $1$.
Call
$$\begin{array}{rcl}
\mathfrak{S}_1 & =& \{Q_2(G_{2j-1},G_{2j}): j \in [1,v+w+1]  \},\\
\mathfrak{S}_2 & =& \{Q_1(H_{2h-1}, H_{2h}): h \in [1,4vw+3(v+w)+1] \}.
\end{array}$$
Then, $\mathfrak{S}=\mathfrak{S}_1\dcup \mathfrak{S}_2$ consists of $4vw+4(v+w)+2$ zero-sum blocks of size $2\times 4$ and
$\E(\mathfrak{S})=\pm (\mathcal{M}_1\dcup \mathcal{M}_3\dcup \mathcal{M}_{5})$.
Hence, we replace the instances of $Z_4$ with the arrays in  $\mathfrak{S}$.
\end{proof}

\begin{lem}\label{a3b3-2}
Suppose that $a,b\geqslant 7$ are such that $a\equiv b\equiv 3 \pmod 4$.
There exists an $\SMAS(a,b;2c)$ for all $c\geqslant 1$ such that $c\equiv 1 \pmod 4$.
\end{lem}

\begin{proof}
By Lemma \ref{77-1} we can write $c=4t+5$ where $t\geqslant 0$. 
Our $\SMAS(a,b;2c)$ will consist of $8t+8$ arrays of type \eqref{97}
and two arrays of type \eqref{97special}.
Write $a=4v+7$ and $b=4w+7$, where $v,w\geqslant 0$.
Hence, we will construct
$8t+8$ blocks $X_3$, 
$4(v+w)(4t+5)+    32t+36 $ blocks $X_2$, 
$2(8t+8)(v+1)(w+1)+2(w+(2v+1)(w+1))$ blocks  $Z_4$,
$2$ blocks $Y_5$ and $2$ blocks $W_2$.
So, take
$$\begin{array}{rcl}
\mathfrak{A} & = & \mathfrak{A}(4t+5,\;  4(v+w)(4t+5)+  32t+40   ,\; 4t+4 ),\\
\mathfrak{B} & = & \mathfrak{B}(9 ,\;  (16vw+ 24(v+w))(4t+5)+164t+201, \;  4(v+w)(4t+5)+    32t+36      ),\\
\mathfrak{C} & =& \mathfrak{C}( 8(v+w)(4t+5)+112t+132   ,\; (16vw+24(v+w))(4t+5)+ 164t+202 , \;  \\
&&8(v+w)(4t+5)+  112t+142 )
\end{array}$$
as in Lemmas \ref{BlockA}, \ref{BlockB} and \ref{BlockC}, respectively.
We have $$\E(\mathfrak{A}) \dcup\E(\mathfrak{B})\dcup \E(\mathfrak{C})=
\pm [1, (16vw+ 28(v+w)) (4t+5)+ 196t+245 ]\setminus \pm \left(\mathcal{M}_1\dcup \ldots \dcup \mathcal{M}_{12} \right),$$
where
$$\begin{array}{rcl}
\mathcal{M}_1 & =& [12, 8t+10]_2, \\
\mathcal{M}_2 & =& [40t+44, 40t+46]_2,\\
\mathcal{M}_3 & =& [40t+50, 56t+62]_4, \\
\mathcal{M}_4 & =& [56t+64,   8(v+w)(4t+5)+ 80t+94]_2, \\
\mathcal{M}_5 & =& [ 8(v+w)(4t+5)+80t+96 ,  8(v+w)(4t+5)+80t+99 ],\\
\mathcal{M}_6 & =& [ 8(v+w)(4t+5)+80t+100,   8(v+w)(4t+5)+112t+130 ]_2,\\
\mathcal{M}_7 & =& [ 8(v+w)(4t+5)+80t+103,   8(v+w)(4t+5)+112t+131]_4,\\
\mathcal{M}_8 & =& [ 8(v+w)(4t+5)+112t+133,   8(v+w)(4t+5)+112t+134],\\
\mathcal{M}_9 & =&  [ 8(v+w)(4t+5)+112t+135,  8(v+w)(4t+5)+ 112t+137]_2,\\
\mathcal{M}_{10} & =&  [ 8(v+w)(4t+5)+112t+139,  8(v+w)(4t+5)+112t+140],\\
\mathcal{M}_{11} & =& [ 8(v+w)(4t+5)+ 112t+146,  8(v+w)(4t+5)+132t+165],\\
\mathcal{M}_{12} & =& [ 8(v+w)(4t+5)+132t+166,  (16vw+20(v+w))(4t+5)+132t+165].
\end{array}$$
So, we replace the $8t+8$ instances of $X_3$ with the arrays in  $\mathfrak{A}$,  
the $2$ instances of $Y_5$ with the arrays in  $\mathfrak{C}$ and the $ 4(v+w)(4t+5)+ 32t+36$ instances of $X_2$ with the arrays in $\mathfrak{B}$.

Define 
$$\mathfrak{R} = \{U_2(\mathcal{M}_8), U_2(\mathcal{M}_{10}) \}.$$
Then, $\E(\mathfrak{R})=\pm(\mathcal{M}_8\dcup \mathcal{M}_{10})$: we replace the $2$ instances of $W_2$ with the arrays
in $\mathfrak{R}$.

The set $\mathcal{M}_3 \dcup\mathcal{M}_7$  can be written as a disjoint union 
$F_1\dcup \ldots \dcup F_{2\alpha}$, where each $\alpha=3t+3$ and $F_i$ is
a $2$-set of type $4$.
The set $\mathcal{M}_1\dcup \mathcal{M}_2 \dcup \mathcal{M}_{4}\dcup \mathcal{M}_{6}\dcup \mathcal{M}_9$
can be written as a disjoint union  $G_1\dcup \ldots \dcup G_{2\beta}$, where $\beta= (v+w)(4t+5)+  8t+9$ and each $G_j$ is
a $2$-set of type $2$.
Finally, the set $\mathcal{M}_5\dcup \mathcal{M}_{11}\dcup \mathcal{M}_{12}$ can be written as a disjoint union 
$H_1\dcup \ldots \dcup H_{2\gamma}$, where $\gamma= (4vw+3(v+w))(4t+5)+ 5t+6$ and each $H_h$ is a $2$-set of type $1$.
Call
$$\begin{array}{rcl}
\mathfrak{S}_1 & =& \left\{Q_3(F_{2i-1},F_{2i} ): i \in [1,3t+3] \right\},\\
\mathfrak{S}_2 & =& \left\{Q_2(G_{2j-1},G_{2j} ): j \in [1,(v+w)(4t+5)+  8t+9] \right\},\\
\mathfrak{S}_3 & =& \left\{Q_1(H_{2h-1},H_{2h} ): h \in [1,(4vw+3(v+w))(4t+5)+ 5t+6] \right\}.
\end{array}$$
Then $\mathfrak{S}=\mathfrak{S}_1\dcup \mathfrak{S}_2\dcup \mathfrak{S}_3$ consists of $(4vw+4(v+w))(4t+5)+ 16t+18$ zero-sum blocks of size $2\times 4$ and
$\E(\mathfrak{S})=\pm (\mathcal{M}_1\dcup \ldots\dcup \mathcal{M}_{7}\dcup \mathcal{M}_9\dcup \mathcal{M}_{11}\dcup \mathcal{M}_{12})$.
Hence, we replace the instances of $Z_4$ with the arrays in  $\mathfrak{S}$.
\end{proof}

\subsection{Subcase \texorpdfstring{$c\equiv 2 \pmod 4$}{c ≡ 2 (mod 4)}}

\begin{lem}\label{59-2n}
Suppose that $a\geqslant 5$ and $b\geqslant 9$ are such that  
$a\equiv  1\pmod 4$ and $b$ is odd.
There exists an $\SMAS(a,b;2c)$ for all $c\geqslant 2$ such that $c\equiv 2 \pmod 4$.
\end{lem}

\begin{proof}
Write $c=4t+2$ where $t\geqslant 0$.
Also, write $a=4v+5$ and $b=4w+9+2\delta$, where $v,w\geqslant 0$ and $\delta \in \{0,1\}$. 
If $b\equiv 1\pmod 4$, then our $\SMAS(a,b;2c)$ will consist of $8t+2$ arrays of type \eqref{97b}
and two arrays of type  \eqref{extra9}; 
if $b\equiv 3 \pmod 4$, it will consist of $8t+2$ arrays of type \eqref{97}
and two arrays of type  \eqref{97special}.
So, we will construct
$8t+2$ blocks $X_3$, 
$(8(v+w)+4\delta)(2t+1) + 32t+12 $ blocks $X_2$, 
$(4t+2)(v(6+2\delta)+w(4v+2))+(16t+6)\delta$ blocks  $Z_4$,
and $2$ blocks $Y_5$.
Furthermore,  we will also construct
$8t+2$ blocks $Z_6$ and  $2$ blocks $W_4$ if $\delta=0$,  or
$2$ blocks $W_2$ if $\delta=1$.
So, take
$$\begin{array}{rcl}
 \mathfrak{A} & =& \mathfrak{A}(4t+5 ,\;
 (8(v+w)+4\delta)(2t+1) +32t+16  ,\; 
 4t+1),\\
 \mathfrak{B} & =&  \mathfrak{B}(9 ,\;
 (32v(w+1)+32(v+w)+16\delta (v+1))(2t+1)+ 148t+70,\;\\
 &&(8(v+w)+4\delta)(2t+1) + 32t+12  ),\\
 \mathfrak{C} & =& \mathfrak{C}(
 (16(v+w)+8\delta)(2t+1) +80t+37, \;\\
 &&(32v(w+1)+32(v+w)+16 \delta (v+1))(2t+1)+  148t+71, \;\\
 &&   (16(v+w) +8\delta)(2t+1)+      116t+55  )
 \end{array}$$
 as in Lemmas \ref{BlockA}, \ref{BlockB} and \ref{BlockC}, respectively.
 We have $$\begin{array}{rcl}\E(\mathfrak{A}) \dcup\E(\mathfrak{B})\dcup \E(\mathfrak{C}) & =& \pm [1, (32v(w+1)+40(v+w)+  4\delta (4v+5) )(2t+1)+180t+90]\setminus \\
 &&\pm \left(\mathcal{M}_1\dcup \ldots \dcup \mathcal{M}_{12} \right),
 \end{array}$$
 where
$$\begin{array}{rcl}
\mathcal{M}_1 & =& [12, 8t+10]_2,\\
\mathcal{M}_2 & =& [40t+20,40t+26]_2,\\
\mathcal{M}_3 & =& [40t+28,  56t+32]_4,\\
\mathcal{M}_4 & =& [56t+34, 8(2v+2w+\delta)(2t+1)+ 80t+36]_2,\\
\mathcal{M}_5 & =&[8(2v+2w+\delta)(2t+1)+80t+38, 
8(2v+2w+\delta)(2t+1)+80t+39],\\
\mathcal{M}_6 & =&[8(2v+2w+\delta)(2t+1)+80t+40,
8(2v+2w+\delta)(2t+1)+80t+42]_2,\\
\mathcal{M}_7 & =&[8(2v+2w+\delta)(2t+1)+80t+44,
8(2v+2w+\delta)(2t+1)+80t+45],\\
\mathcal{M}_8 & =&[8(2v+2w+\delta)(2t+1)+80t+48,
8(2v+2w+\delta)(2t+1)+80t+49],\\
\mathcal{M}_9 & =&[8(2v+2w+\delta)(2t+1)+80t+50,
8(2v+2w+\delta)(2t+1)+112t+52]_2,\\
\mathcal{M}_{10} & =&[8(2v+2w+\delta)(2t+1)+80t+53,
8(2v+2w+\delta)(2t+1)+112t+49]_4,\\
\mathcal{M}_{11} & =&[ 8(2v+2w+\delta)(2t+1)+112t+53,
8(2v+2w+\delta)(2t+1)+116t+54],\\
\end{array}$$

$$\begin{array}{rcl}
\mathcal{M}_{12} & =&[8(2v+2w+\delta)(2t+1)+116t+59,\\
&&
(12(2v+2w+\delta)+16v(2w+ 2+\delta))(2t+1)+ 116t+58].
\end{array}$$
We replace the $8t+2$ instances of $X_3$ with the arrays in  
$\mathfrak{A}$,  
the $2$ instances of $Y_5$ with the arrays in  $\mathfrak{C}$ and 
the $(8(v+w)+4\delta)(2t+1) + 32t+12$ instances of $X_2$ with the arrays in 
$\mathfrak{B}$.

Write
$$  \mathcal{M}_{9}=\mathcal{M}_{9}^1\dcup  \mathcal{M}_{9}^2
 \dcup \mathcal{M}_9^3\equad
 \mathcal{M}_{11}=\mathcal{M}_{11}^1\dcup  \mathcal{M}_{11}^2
 \dcup \mathcal{M}_{11}^3, $$
where
 $$\begin{array}{rcl}
\mathcal{M}_{9}^1 & =& [8(2v+2w+\delta)(2t+1)+80t+50,
8(2v+2w+\delta)(2t+1)+104t+52]_2,\\
\mathcal{M}_{9}^2 & =& [8(2v+2w+\delta)(2t+1)+104t+54,
8(2v+2w+\delta)(2t+1)+112t+50]_4,\\
\mathcal{M}_{9}^3 & =& [8(2v+2w+\delta)(2t+1)+104t+56,
8(2v+2w+\delta)(2t+1)+112t+52]_4,\\
\mathcal{M}_{11}^1 & =&[ 8(2v+2w+\delta)(2t+1)+112t+53,
8(2v+2w+\delta)(2t+1)+112t+54],\\ 
\mathcal{M}_{11}^2 & =& [ 8(2v+2w+\delta)(2t+1)+112t+55,
8(2v+2w+\delta)(2t+1)+116t+53]_2,\\
\mathcal{M}_{11}^3 & =& [ 8(2v+2w+\delta)(2t+1)+112t+56,
8(2v+2w+\delta)(2t+1)+116t+54]_2.
\end{array}$$
The set $\mathcal{M}_3\dcup \mathcal{M}_9^2 \dcup\mathcal{M}_9^3\dcup \mathcal{M}_{10}$  can be written as a disjoint union 
$F_1\dcup \ldots \dcup F_{\alpha}$, where $\alpha=8t+1$  and each $F_i$ is a $2$-set of type $4$.
The set $\mathcal{M}_1\dcup\mathcal{M}_2\dcup \mathcal{M}_4\dcup \mathcal{M}_{9}^1 \dcup  \mathcal{M}_{11}^2\dcup \mathcal{M}_{11}^3$
can be written as a disjoint union 
$G_1\dcup \ldots \dcup G_{2\beta+2(8t+2)}$, where $\beta=
(4t+2) (v+w)+(2t+1)\delta$ and each $G_j$ is
a $2$-set of type $2$.
Finally, the set $\mathcal{M}_7\dcup  \mathcal{M}_8\dcup \mathcal{M}_{11}^1\dcup \mathcal{M}_{12}$ can be written as a disjoint union
$H_1\dcup\ldots\dcup H_{2\gamma+3}$, where $\gamma=
(4t+2) (v (5+2\delta)+ w(4v+1)) + (2t+1)\delta$ and each $H_h$ is a $2$-set of type $1$.

Call
$$\begin{array}{rcl}
\mathfrak{Q}_1 & =& \left\{Q_2(G_{2j-1},G_{2j} ): j \in [1, \beta] \right\},\\
\mathfrak{Q}_2 & =& \left\{Q_1(H_{2h-1},H_{2h} ): h \in [1,\gamma]  \right\}.
 \end{array}$$
Then $\mathfrak{Q}=\mathfrak{Q}_1\dcup \mathfrak{Q}_2$ consists of $\beta+\gamma=(4t+2) (v (6+2\delta)+ w(4v+2)) + (4t+2)\delta$ zero-sum blocks of size $2\times 4$.

Suppose $\delta=0$. 
First of all, we replace the instances of $Z_4$ with the arrays in  
$\mathfrak{Q}$. Then, we need still to construct $8t+2$ blocks $Z_6$ and $2$ blocks $W_4$. To this purpose, define 
$$\begin{array}{rcl}
\mathfrak{R}_1 & =& \left\{R_2(F_i,G_{2\beta+2i-1},G_{2\beta+2i} ): i \in [1,8t+1] \right\},\\
\mathfrak{R}_2 & =& \left\{R_1(G_{2\beta+2(8t+1)+1}, H_{2\gamma+1},H_{2\gamma+2} )  \right\},\\
\mathfrak{S} & = & \{U_1(\mathcal{M}_6, \mathcal{M}_5), U_1(G_{2\beta+2(8t+2)},H_{2\gamma+3} ) \}.
\end{array}$$
Then  $\mathfrak{R}=\mathfrak{R}_1\dcup \mathfrak{R}_2$ consists
of $8t+2$ zero-sum blocks of size $2\times 6$.
Hence, we replace the instances of $Z_6$ with the arrays in  
$\mathfrak{R}$ and the instances of $W_4$ with the arrays in  
$\mathfrak{S}$.

Next, suppose $\delta=1$. Write $\mathcal{M}_5 \dcup \mathcal{M}_6$
as $\mathcal{N}_1\dcup \mathcal{N}_2$, where
$$\begin{array}{rcl}
\mathcal{N}_1 & =&[8(2v+2w+\delta)(2t+1)+80t+38, 
8(2v+2w+\delta)(2t+1)+80t+42]_4,\\
\mathcal{N}_2 & =&[8(2v+2w+\delta)(2t+1)+80t+39,
8(2v+2w+\delta)(2t+1)+80t+40].
\end{array}$$
Define
$$\begin{array}{rcl}
\mathfrak{Q}_3 & =& \left\{Q_3(F_{2i-1},F_{2i} ): j \in [1, 4t] \right\}\dcup\{Q_3(F_{8t+1},\mathcal{N}_1)\},\\
\mathfrak{Q}_4 & =& \left\{Q_2(G_{2\beta+2j-1},G_{2\beta+2j} ): i \in [1, 8t+2] \right\},\\
\mathfrak{Q}_5 & =& \left\{Q_1(H_{2\gamma+1},H_{2\gamma+2} ) \right\}.
 \end{array}$$
Then $\mathfrak{Q}'=\mathfrak{Q}_3\dcup \mathfrak{Q}_4\dcup\mathfrak{Q}_5$ consists of $12t+4$ zero-sum blocks of size $2\times 4$.
So, we replace the instances of $Z_4$ with the arrays in  
$\mathfrak{Q}\dcup \mathfrak{Q}'$.  To conclude,  we construct the two  blocks $W_2$ by taking $U_2(\mathcal{N}_2)$ and $U_2(H_{2\gamma+3})$.
\end{proof}

\begin{lem}\label{97-c2}
Suppose that $a\geqslant 9$ is such that  $a\equiv 1 \pmod 4$.
There exists an $\SMAS(a,7;2c)$ for all $c\geqslant 2$ such that 
$c\equiv 2 \pmod 4$.
\end{lem}

\begin{proof}
Write $c=4t+2$ where $t\geqslant 0$.
Our $\SMAS(a,7;2c)$ will consist of $8t+2$ arrays of type \eqref{97sm} and two arrays of type  \eqref{97supsp}.
Write $a=4v+9$, where $v\geqslant 0$.
Hence, we will construct
$8t+2$ blocks $X_3$, 
$ 8v(2t+1) +  40t+16  $ blocks $X_2$,
  $(8t+2)(2v+3)+2(2v+2)$  blocks  $Z_4$,
$2$ blocks $Y_5$ and $2$ blocks $W_2$.
So, take
 $$\begin{array}{rcl}
\mathfrak{A} & =& \mathfrak{A}(4t+5 ,\; 8v(2t+1)+ 40t+20  ,\; 4t+1),\\
 \mathfrak{B} & =&   \mathfrak{B}(9 ,\; 48v(2t+1)+ 212t+102 ,\; 8v(2t+1) +  40t+16  ),\\
 \mathfrak{C} & = & \mathfrak{C}(16v(2t+1)+   96t+42 ,\; 48v(2t+1)+ 212t+103 ,\; 16v(2t+1)+  128t+61 )
 \end{array}$$
 as in Lemmas \ref{BlockA},  \ref{BlockB} and  \ref{BlockC}, respectively.
 We have $$\E(\mathfrak{A}) \dcup\E(\mathfrak{B})\dcup \E(\mathfrak{C})=\pm [1, 56v(2t+1)+252t+126]
 \setminus \pm \left(\mathcal{M}_1\dcup \ldots \dcup  \mathcal{M}_{11} \right),$$
where
$$\begin{array}{rcl}
\mathcal{M}_1 & =& [12,8t+10]_2, \\
\mathcal{M}_2 & =& [40t+20, 40t+26]_2, \\
\mathcal{M}_3 & =& [40t+28, 56t+32]_4,  \\
\mathcal{M}_4 & =& [56t+34,  16v(2t+1)+96t+40]_2 ,\\
\mathcal{M}_5 & =& [16v(2t+1)+96t+44, 16v(2t+1)+96t+45], \\
\mathcal{M}_6 & =& [16v(2t+1)+96t+47,  16v(2t+1)+96t+49]_2, \\
\mathcal{M}_7 & =& [ 16v(2t+1)+96t+50,  16v(2t+1)+96t+52]_2, \\
\mathcal{M}_8 & =& [ 16v(2t+1)+96t+53,  16v(2t+1)+ 128t+57]_4, \\
\mathcal{M}_9 & =& [16v(2t+1)+96t+54, 16v(2t+1)+128t+60]_2,\\
\mathcal{M}_{10} & =& [16v(2t+1)+128t+65, 16v(2t+1)+172t+86],\\
\mathcal{M}_{11} & =& [16v(2t+1)+172t+87,  40v(2t+1)+172t+86].
\end{array}$$
So, we replace the $8t+2$ instances of $X_3$ with the arrays in  $\mathfrak{A}$,  
the $2$  instances of $Y_5$ with the arrays in  $\mathfrak{C}$ and the $8v(2t+1) +  40t+16$ instances of $X_2$ with the arrays in $\mathfrak{B}$.

Write
 $$\mathcal{M}_{10}=\mathcal{M}_{10}^1\dcup  \mathcal{M}_{10}^2,$$
where
 $$\begin{array}{rcl}
\mathcal{M}_{10}^1 & =& [16v(2t+1)+128t+65, 16v(2t+1)+128t+66 ],\\
\mathcal{M}_{10}^2 & =& [16v(2t+1)+128t+67, 16v(2t+1)+172t+86].
\end{array}$$
Define 
$$\mathfrak{R}=\{U_2(\mathcal{M}_{5}), U_2(\mathcal{M}_{10}^1) \}.$$
 Then, $\E(\mathfrak{R})=\pm(\mathcal{M}_5\dcup \mathcal{M}_{10}^1)$: we replace the $2$ instances of $W_2$ with the arrays
 in $\mathfrak{R}$. 
 
The set $\mathcal{M}_3\dcup \mathcal{M}_8$  can be written as a disjoint union 
$F_1\dcup \ldots \dcup F_{2\alpha}$, where $\alpha=3t+1$  and each $F_i$ is a $2$-set of type $4$.
The set $\mathcal{M}_1\dcup \mathcal{M}_2\dcup \mathcal{M}_{4}\dcup \mathcal{M}_6\dcup  \mathcal{M}_{7}\dcup \mathcal{M}_{9}$
can be written as a disjoint union 
$G_1\dcup \ldots \dcup G_{2\beta}$, where $\beta= 2v(2t+1) +  10t+4 $ and each $G_j$ is
a $2$-set of type $2$.
Finally, the set $\mathcal{M}_{10}^2\dcup \mathcal{M}_{11}$ can be written as a disjoint union
$H_1\dcup \ldots \dcup H_{2\gamma}$, where $\gamma= 6v(2t+1)+ 11t+5$ and each $H_h$ is a $2$-set of type $1$.
Call
$$\begin{array}{rcl}
\mathfrak{S}_1 & =& \left\{Q_3(F_{2i-1},F_{2i} ): i \in [1, 3t+1] \right\},\\
\mathfrak{S}_2 & =& \left\{Q_2(G_{2j-1},G_{2j} ): j \in [1, 2v(2t+1) +  10t+4] \right\},\\
\mathfrak{S}_3 & =& \left\{Q_1(H_{2h-1},H_{2h} ): h \in [1, 6v(2t+1)+ 11t+5]  \right\}.
 \end{array}$$
Then $\mathfrak{S}=\mathfrak{S}_1\dcup \mathfrak{S}_2\dcup \mathfrak{S}_3$ consists of $8v(2t+1)+24t+10$ zero-sum blocks of size $2\times 4$ and
$\E(\mathfrak{S})=\pm  (\mathcal{M}_1\dcup \ldots\dcup\mathcal{M}_4 \dcup \mathcal{M}_6\dcup\ldots\dcup \mathcal{M}_{9}\dcup \mathcal{M}_{10}^2\dcup 
\mathcal{M}_{11})$.
 Hence, we replace the instances of $Z_4$ with the arrays in  $\mathfrak{S}$. 
\end{proof}

\begin{lem}\label{ab3-2}
Suppose that $a,b \geqslant 7$ are such that $a\equiv b\equiv 3 \pmod 4$.
There exists an $\SMAS(a,b;2c)$ for all $c\geqslant 2$ such that $c\equiv 2 \pmod 4$.
\end{lem}

\begin{proof}
Write $c=4t+2$ where $t\geqslant 0$. Our $\SMAS(a,b;2c)$ will consist
of $8t+2$ arrays of type \eqref{97} and two arrays of type
\eqref{97special}.
Write $a=4v+7$ and $b=4w+7$, where $v,w\geqslant 0$.
Hence, we will construct
$8t+2$ blocks $X_3$, 
$8(v+w)(2t+1)+  32t+12 $ blocks $X_2$,
 $ 2(8t+2)(v+1)(w+1)+ 2(w+(2v+1)(w+1))$    blocks  $Z_4$,
$2$ blocks $Y_5$ and $2$ blocks $W_2$.
So, take
$$\begin{array}{rcl}
\mathfrak{A} &  = & \mathfrak{A}(4t+5 ,\; 8(v+w)(2t+1)+ 32t+16  ,\; 4t+1),\\
 \mathfrak{B} & =&  \mathfrak{B}(9 ,\; (32vw+48(v+w))(2t+1)+ 164t+78
 ,\;  8(v+w)(2t+1)+  32t+12  ),\\
 \mathfrak{C} & =& \mathfrak{C}( 16(v+w)(2t+1)+ 80t+34 ,\; (32vw+48(v+w))(2t+1)+ 164t+79 , \;   \\
 && 16(v+w)(2t+1)+   112t+55)
 \end{array}$$
 as in Lemmas \ref{BlockA}, \ref{BlockB} and  \ref{BlockC}, respectively.
 We have $$\E(\mathfrak{A}) \dcup\E(\mathfrak{B})\dcup \E(\mathfrak{C})=\pm [1, (32vw+56(v+w))(2t+1)+196t+98 ]
 \setminus \pm \left(\mathcal{M}_1\dcup \ldots \dcup \mathcal{M}_{12} \right),$$
 where
 $$\begin{array}{rclcrcl}
\mathcal{M}_1 & =& [16(v+w)(2t+1)+ 80t+36, 16(v+w)(2t+1)+ 80t+37] , \\
\mathcal{M}_2 & =& [16(v+w)(2t+1)+ 80t+45, 16(v+w)(2t+1)+ 80t+46], \\
\mathcal{M}_3 & =&  [12,8t+10]_2, \\
\mathcal{M}_4 & =& [40t+20,  40t+26]_2,  \\
\mathcal{M}_5 & =& [40t+28, 56t+32]_4,\\
\mathcal{M}_6 & =& [56t+34, 16(v+w)(2t+1)+  80t+32]_2, \\
\mathcal{M}_7 & =& [16(v+w)(2t+1)+ 80t+39, 16(v+w)(2t+1)+ 80t+41]_2,\\
\mathcal{M}_8 & =& [16(v+w)(2t+1)+ 80t+42, 16(v+w)(2t+1)+ 80t+44]_2, \\
\mathcal{M}_9 & =&[16(v+w)(2t+1)+ 80t+48, 16(v+w)(2t+1)+ 112t+54]_2,\\
\mathcal{M}_{10} & =& [16(v+w)(2t+1)+ 80t+49, 16(v+w)(2t+1)+ 112t+53 ]_4, \\
\mathcal{M}_{11} & =& [16(v+w)(2t+1) +112t+59, 16(v+w)(2t+1) +132t+66],\\
\mathcal{M}_{12} & =& [16(v+w)(2t+1)+132t+67,   (32vw +40(v+w))(2t+1)+ 132t+66].
\end{array}$$
So, we replace the $8t+2$ instances of $X_3$ with the arrays in  $\mathfrak{A}$,  
the $2$  instances of $Y_5$ with the arrays in  $\mathfrak{C}$ and the $ 8(v+w)(2t+1)+  32t+12 $ instances of $X_2$ with the arrays in $\mathfrak{B}$.

Define 
$$\mathfrak{R}=\{U_2(\mathcal{M}_1), U_2(\mathcal{M}_2) \}.$$
 Then, $\E(\mathfrak{R})=\pm(\mathcal{M}_1\dcup \mathcal{M}_{2})$:  we replace the $2$ instances of $W_2$ with the arrays in $\mathfrak{R}$. 
 
The set $\mathcal{M}_5\dcup \mathcal{M}_{10}$  can be written as a disjoint union 
$F_1\dcup \ldots \dcup F_{2\alpha}$, where $\alpha=3t+1$  and each $F_i$ is a $2$-set of type $4$.
The set $\mathcal{M}_3\dcup \mathcal{M}_4\dcup \mathcal{M}_{6}\dcup \mathcal{M}_7\dcup  \mathcal{M}_{8}\dcup \mathcal{M}_{9}$
can be written as a disjoint union 
$G_1\dcup \ldots \dcup G_{2\beta}$, where $\beta= 2(v+w)(2t+1)+   8t+3 $ and each $G_j$ is
a $2$-set of type $2$.
Finally, the set $\mathcal{M}_{11}\dcup \mathcal{M}_{12}$ can be written as a disjoint union
$H_1\dcup \ldots \dcup H_{2\gamma}$, where $\gamma= 
(8vw+ 6(v+w)) (2t+1)+ 5t+2$ and each $H_h$ is a $2$-set of type $1$.
Call
$$\begin{array}{rcl}
\mathfrak{S}_1 & =& \left\{Q_3(F_{2i-1},F_{2i} ): i \in [1,3t+1] \right\},\\
\mathfrak{S}_2 & =& \left\{Q_2(G_{2j-1},G_{2j} ): j \in [1,2(v+w)(2t+1)+8t+3] \right\},\\
\mathfrak{S}_3 & =& \left\{Q_1(H_{2h-1},H_{2h} ): h \in [1,(8vw+ 6(v+w)) (2t+1)+ 5t+2]  \right\}.
 \end{array}$$
Then $\mathfrak{S}=\mathfrak{S}_1\dcup \mathfrak{S}_2\dcup \mathfrak{S}_3$ consists of $(8vw+8(v+w))(2t+1)+  16t+6$ zero-sum blocks of size $2\times 6$ and
$\E(\mathfrak{S})=\pm  (\mathcal{M}_3\dcup \ldots\dcup \mathcal{M}_{12})$.
 Hence, we replace the instances of $Z_4$ with the arrays in  $\mathfrak{S}$.
\end{proof}

\section{The case \texorpdfstring{$\{ a,b \}=\{5,7\}$}{{a, b}={5, 7}} and the case 
\texorpdfstring{$c\equiv 0\pmod 4$}{c ≡ 0 (mod 4)}}\label{sm}

In this section, we consider the case of an $\SMA(a,b;2c)$ where $c\equiv 0 \pmod 4$ and  the exceptional case of an $\SMAS(5,7;2c)$ that require slightly different constructions.

\begin{lem}\label{BlockA-2}
Given four nonnegative integers $\alpha,\beta,\gamma,u$ such that
$\alpha+1\geqslant 2u$ and $2\gamma \geqslant \alpha+8u$,
there exists a set $\overline{\mathfrak{A}}=\overline{\mathfrak{A}}(\alpha,\beta,\gamma,u)$ consisting of $4u$ square zero-sum blocks of size $3$ such that
$\E(\overline{\mathfrak{A}})$ consists of the following disjoint subsets:
$$\begin{array}{rcl}
\E(\overline{\mathfrak{A}})  & = & \pm[2\alpha+2,  2\alpha+   16u]_2 \dcup  \pm[4\alpha+12u+4, 4\alpha+20u]_4\dcup\\
&&  \pm[2\beta+1,  2\beta+2u-1 ]_2 \dcup  \pm[2\beta+4u+1, 2\beta+6u-1]_2 \dcup  \\
&& \pm [2\alpha+2\beta+4u+3, 2\alpha+2\beta+8u-1 ]_4 \dcup \\
&& \pm [2\alpha+2\beta+12u+3, 2\alpha+2\beta+16u-1]_4\dcup 
 \pm[4\gamma+2,  4\gamma+8u-2]_4\dcup \\
&&  \pm[2\alpha+4\gamma+8u+2,   2\alpha+4\gamma+10u]_2 \dcup 
 \pm[2\alpha+4\gamma+16u+2,   2\alpha+4\gamma+18u]_2.
\end{array}$$
\end{lem}

\begin{proof}
For all $i\in [0,u-1]$, define
$$\begin{array}{rcl}
A_{2i}& =& \begin{array}{|c|c|c|}\hline
4\alpha+12u+4+4i    & -(2\alpha+4u+2+2i)       & -(2\alpha+8u+2+2i) \\ \hline
-(2\alpha+2+2i)     &  2\alpha+2\beta+4u+3+4i  & -(2\beta+4u+1+2i)\\ \hline
-(2\alpha+12u+2+2i) & -(2\beta+1+2i)           &  2\alpha+2\beta+12u+3+4i \\\hline
\end{array},\\\\[-8pt]
A_{2i+1}& =& \begin{array}{|c|c|c|}\hline
4\alpha+16u+4+4i    & -(2\alpha+6u+2+2i)        & -(2\alpha+10u+2+2i) \\\hline
-(2\alpha+2u+2+2i)  &  2\alpha+4\gamma+10u-2i & -(4\gamma+8u-2-4i)\\\hline
-(2\alpha+14u+2+2i) & -(4\gamma+4u-2-4i)        &  2\alpha+4\gamma+18u-2i\\\hline
\end{array}.
\end{array}$$
The set $\overline{\mathfrak{A}}=\{+A_{2i},-A_{2i}, +A_{2i+1}, -A_{2i+1} :  i \in [0,u-1]\}$ has the required properties.
\end{proof}

\begin{lem}\label{57-0}
Suppose that $b\geqslant 7$ is an odd integer.
There exists an $\SMAS(5,b;2c)$ for all $c\geqslant 4$ such that $c\equiv 0 \pmod 4$.
\end{lem}

\begin{proof}
Write $c=4t+4$ where $t\geqslant 0$.
Our $\SMAS(5,b;2c)$ will consist
of $8t+8$ arrays of type \eqref{59}
if $b\equiv 1 \pmod 4$,
or of type \eqref{57}
if $b\equiv 3 \pmod 4$, 
where each $X_r$ is a zero-sum block of size $r\times 3$ and 
each $Z_\ell$ is a zero-sum block of size $2\times \ell$.

Write $b=4w+7+2\delta$, where $w\geqslant 0$ and 
$\delta\in\{0,1\}$.
Hence, we will construct
$8t+8$ blocks $X_3$, 
$8(2w+3+\delta)(t+1)$
blocks $X_2$,
$8(w+1-\delta)(t+1)$ blocks  $Z_4$
and $8\delta(t+1)$  blocks $Z_6$.
Take
$$\begin{array}{rcl}
\overline{\mathfrak{A}} & = & 
\overline{\mathfrak{A}}(4(t+1),\; 
8(2w+3+\delta)(t+1),\;  10(t+1)  ,\; 2(t+1)),\\
\mathfrak{B} & = & \mathfrak{B}(1,\; 
4(16w+29+8\delta)(t+1), \;  8(2w+3+\delta)(t+1) )   
  \end{array}$$
as in Lemmas \ref{BlockA-2} and \ref{BlockB}, respectively.
Then, 
$$\E(\overline{\mathfrak{A}})\dcup \E(\mathfrak{B})=\pm [1, 
40(2w+\delta)(t+1)+140t+140]
\setminus \pm \left(\mathcal{M}_1\dcup \ldots \dcup \mathcal{M}_{10} \right),$$
where
$$\begin{array}{rcl}
\mathcal{M}_1 & = & [2,8t+8]_2,\\
\mathcal{M}_2 & = & [56t+58, 64t+64]_2,\\
\mathcal{M}_3 & = & [68t+70,  80t+80]_2,\\
\mathcal{M}_4 & = & [84t+86,  16(2w+\delta)(t+1)+ 88t+88 ]_2,\\
\mathcal{M}_5 & = & [ 16(2w+\delta)(t+1) +52t+53,  16(2w+\delta)(t+1) +56t+55]_2,\\
\mathcal{M}_6 & = &  [16(2w+\delta)(t+1) +60t+61, 16(2w+\delta)(t+1)+64t+63]_2,\\
\mathcal{M}_7 & = & [16(2w+\delta)(t+1)+64t+65, 
16(2w+\delta)(t+1)+72t+69]_4,\\
\mathcal{M}_8 & = & [16(2w+\delta)(t+1) +72t+73, 16(2w+\delta)(t+1) +80t+79]_2,\\
\mathcal{M}_9 & = & [16(2w+\delta)(t+1)+80t+81, 16(2w+\delta)(t+1)+88t+85]_4,\\
\mathcal{M}_{10} & = & [ 16(2w+\delta)(t+1) +88t+89, 24(2w+\delta)(t+1) +92t+92].
\end{array}$$
So, we replace the $8t+8$ instances of $X_3$ with the arrays in  
$\overline{\mathfrak{A}}$ and the $8(2w+3+\delta)(t+1)$ instances of $X_2$ with the arrays in $\mathfrak{B}$.

Next, write
$$\mathcal{M}_1=\mathcal{M}_1^1\dcup \mathcal{M}_1^2,\quad
\mathcal{M}_2=\mathcal{M}_2^1\dcup \mathcal{M}_2^2
\equad 
\mathcal{M}_{10}=\mathcal{M}_{10}^1\dcup \mathcal{M}_{10}^2\dcup \mathcal{M}_{10}^3,
$$
where
$$\begin{array}{rcl}
\mathcal{M}_1^1 & =& [2,8t+6]_4, \\
\mathcal{M}_1^2 & =& [4,8t+8]_4, \\
\mathcal{M}_2^1 & =& [56t+58, 64t+62]_4, \\
\mathcal{M}_2^2 & =& [56t+60, 64t+64]_4, \\
\mathcal{M}_{10}^1 & =&[16(2w+\delta)(t+1) +88t+89,
16(2w+\delta)(t+1) +92t+91]_2,\\
\mathcal{M}_{10}^2 & =&[16(2w+\delta)(t+1) +88t+90,
16(2w+\delta)(t+1) +92t+92]_2,\\
 \mathcal{M}_{10}^3 & =& [
 16(2w+\delta)(t+1) +92t+93, 24(2w+\delta)(t+1) +92t+92]
\end{array}$$

The set $\mathcal{M}_1^1\dcup \mathcal{M}_1^2\dcup \mathcal{M}_2^1\dcup \mathcal{M}_2^2\dcup 
\mathcal{M}_7\dcup \mathcal{M}_9$  can be written as a disjoint union 
$F_1\dcup \ldots \dcup F_{2\alpha}$, where
$\alpha=3t+3$ and each $F_i$ is
a $2$-set of type $4$.
The set $\mathcal{M}_3\dcup \mathcal{M}_4\dcup \mathcal{M}_5 \dcup \mathcal{M}_6\dcup \mathcal{M}_8\dcup \mathcal{M}_{10}^1\dcup \mathcal{M}_{10}^2$
can be written as a disjoint union $G_1\dcup \ldots \dcup G_{2\beta}$,  where
$\beta=(4w+5+2\delta)(t+1)$ and each $G_j$ is a $2$-set of type $2$.
Finally, the set $\mathcal{M}_{10}^3$
can be written as a disjoint union 
$H_1\dcup \ldots \dcup H_{2\gamma}$, where $\gamma=(4w+2\delta)(t+1)$ and each $H_h$ is
a $2$-set of type~$1$.

Call
$$\begin{array}{rcl}
\mathfrak{Q}_1 & =& \left\{Q_2(G_{2j-1},G_{2j} ):
 j \in [1,4w(t+1)] \right\},\\
\mathfrak{Q}_2 & =& \left\{Q_1(H_{2h-1},H_{2h} ):
 h \in [1,4w(t+1)] \right\}.
\end{array}$$
If $\delta=0$, define $\mathfrak{R}=\varnothing$ and 
$\mathfrak{Q}=\mathfrak{Q}_1\dcup \mathfrak{Q}_2\dcup 
\mathfrak{Q}_3\dcup \mathfrak{Q}_4$, where
$$\begin{array}{rcl}
\mathfrak{Q}_3 & =& \left\{Q_3(F_{2i-1},F_{2i} ):
 i \in [1,3t+3] \right\},\\
\mathfrak{Q}_4 & =& \left\{Q_2(G_{2j-1},G_{2j} ):
 j \in [4w(t+1)+1,  (4w+5)(t+1)] \right\}.
 \end{array}$$
If $\delta=1$, define  $\mathfrak{Q}=\mathfrak{Q}_1 \dcup \mathfrak{Q}_2$ and 
$\mathfrak{R}=\mathfrak{R}_1\dcup \mathfrak{R}_2$, where
$$\begin{array}{rcl}
\mathfrak{R}_1 & =& \left\{R_2(F_i, G_{8w(t+1)+2i-1}, G_{8w(t+1)+2i} ):
 i \in [1,6t+6] \right\},\\
\mathfrak{R}_2 & =& \left\{
R_1(G_{(8w+12)(t+1)+j},H_{8w(t+1)+2j-1},H_{8w(t+1)+2j} ):
 j \in [1,2t+2] \right\}.
 \end{array}$$

Then $\mathfrak{Q}$ consists of 
$8(w+1-\delta)(t+1)$ zero-sum blocks of size $2\times 4$, $\mathfrak{R}$ consists of 
$8\delta(t+1)$ zero-sum blocks of size $2\times 6$
and 
$\E(\mathfrak{Q})\dcup \E(\mathfrak{R})=\pm  (\mathcal{M}_1\dcup \ldots\dcup \mathcal{M}_{10})$.
Hence, we replace the  instances of $Z_4$ with the arrays in  $\mathfrak{Q}$
and the  instances of $Z_6$ with the arrays in  $\mathfrak{R}$.
\end{proof}

\begin{lem}\label{57-1}
There exists an $\SMAS(5,7;2c)$ for all $c\geqslant 1$ such that $c\equiv 1 \pmod 4$.
\end{lem}

\begin{proof}
An $\SMAS(5,7;2)$  is shown in Figure \ref{572}. So, write $c=4t+5$ where $t\geqslant 0$.
Our $\SMAS(5,7;2c)$ will consist
of $8t+10$ arrays of type
\begin{equation}\label{570}
\begin{array}{|c|c|c|}\hline
  X_3      &    X_{2}^\T  & X_{2}^\T \\ \hline
  X_{2}    &  \multicolumn{2}{c|}{Z_4} \\ \hline  
\end{array}.
\end{equation}
Hence, we will construct
$8t+10$ blocks $X_3$, 
$6(4t+5)$  blocks $X_2$ and 
$2(4t+5)$   blocks  $Z_4$.

Take
$$\begin{array}{rcl}
\overline{\mathfrak{A}} & = & \overline{\mathfrak{A}}(4t+5 ,\; 6(4t+5) ,\;  10t+12 ,\; 2t+2 ),\\
\mathfrak{B} & = & \mathfrak{B}(1,\;  29(4t+5), \; 6(4t+5)   )   
  \end{array}$$
as in Lemmas \ref{BlockA-2} and \ref{BlockB}, respectively.
Also, define
$$A'=\begin{array}{|c|c|c|}\hline
 10 & -4        & -6 \\\hline
 -2 & 92t+111    & -(92t+109) \\ \hline 
 -8 & -(92t+107) & 92t+115\\ \hline
\end{array}$$
and set $\mathfrak{A}'=\overline{\mathfrak{A}}\dcup \{+A',-A'\}$. 
Then, 
$$\E(\mathfrak{A}')\dcup\E( \mathfrak{B})=\pm [1, 140t+175]
\setminus \pm \left(\mathcal{M}_1\dcup \ldots \dcup \mathcal{M}_{17} \right),$$
 where
 $$\begin{array}{rclcrcl}
\mathcal{M}_1 & =& [12,8t+10]_2, & \quad &
\mathcal{M}_2 & =& [40t+44, 40t+46]_2,\\
\mathcal{M}_3 & =& [52t+65, 56t+63]_2, &&
\mathcal{M}_4 & =& [56t+64, 56t+67],\\
\mathcal{M}_5 & =& [56t+68, 60t+70]_2, &&
\mathcal{M}_6 & =& [60t+72,  64t+75],\\
\mathcal{M}_7 & =& [64t+77, 64t+79]_2, &&
\mathcal{M}_8 & =& [64t+83, 72t+87]_4,\\
\mathcal{M}_9 & =& [ 68t+80, 72t+82]_2, &&
\mathcal{M}_{10} & =& [72t+84, 72t+86]_2,\\
\mathcal{M}_{11} & =& [72t+88, 80t+91], &&
\mathcal{M}_{12} & =& [80t+93, 80t+95]_2,\\
\mathcal{M}_{13} & =& [80t+99, 88t+103]_4, &&
\mathcal{M}_{14} & =& [84t+96, 88t+102]_2,\\
\mathcal{M}_{15} & =& [88t+104, 92t+105], &&
\mathcal{M}_{16} & =& [92t+106, 92t+112]_2,\\
\mathcal{M}_{17} & =& [92t+113, 92t+114].
\end{array}$$
So, we replace the $8t+10$ instances of $X_3$ with the arrays in  $\mathfrak{A}'$ and the $6(4t+5)$ instances of $X_2$ with the arrays in $\mathfrak{B}$.

The set $\mathcal{M}_8\dcup \mathcal{M}_{13}$  can be written as a disjoint union 
$F_1\dcup \ldots \dcup F_{2\alpha}$, where $\alpha=t+1$ and each $F_i$ is a $2$-set of type $4$.
The set $\mathcal{M}_1\dcup \mathcal{M}_2\dcup \mathcal{M}_3\dcup \mathcal{M}_5\dcup\mathcal{M}_7
\dcup \mathcal{M}_{9}\dcup \mathcal{M}_{10}\dcup \mathcal{M}_{12}\dcup \mathcal{M}_{14}\dcup \mathcal{M}_{16}$ 
 can be written as a disjoint union $G_1\dcup \ldots \dcup G_{2\beta}$,  where  $\beta=3t+5$ and each $G_j$ is a $2$-set of type $2$.
Finally, the set $\mathcal{M}_4\dcup \mathcal{M}_6\dcup \mathcal{M}_{11}\dcup \mathcal{M}_{15}\dcup \mathcal{M}_{17}$  can be written as a disjoint union 
 $H_1\dcup \ldots \dcup H_{2\gamma}$, where $\gamma=4t+4$ and each $H_h$ is a $2$-set of type $1$.
Call
 $$\begin{array}{rcl}
\mathfrak{R}_1 & =& \left\{Q_3(F_{2i-1},F_{2i} ):  i \in [1,t+1] \right\},\\
\mathfrak{R}_2 & =& \left\{Q_2(G_{2j-1},G_{2j} ):  j \in [1,3t+5] \right\},\\
\mathfrak{R}_3 & =& \left\{Q_1(H_{2h-1},H_{2h} ):  h \in [1,4t+4] \right\}.
 \end{array}$$
 Then $\mathfrak{R}=\mathfrak{R}_1\dcup \mathfrak{R}_2\dcup  \mathfrak{R}_3$ consists of 
 $8t+10$ zero-sum blocks of size $2\times 4$ and $\E(\mathfrak{R})=\pm  (\mathcal{M}_1\dcup \ldots\dcup \mathcal{M}_{17})$.
 Hence, we replace the instances of $Z_4$ with the arrays in  $\mathfrak{R}$.
\end{proof}

\begin{figure}[ht]
\begin{footnotesize}
$\begin{array}{|c|c|c|c|c|c|c|c|} \hline   
  10 & -8 & -2 & 5 & -5 & 17 & -17 \\ \hline
  -6 & -59 & 65 & -90 & 90 & -96 & 96 \\ \hline
  -4 & 68 & -63 & 84 & -85 & 79 & -79 \\ \hline 
  1 & -89 & 87 & -45 & 47 & 49 & -50 \\ \hline 
  -1 & 88 & -87 & 46 & -47 & -49 & 50 \\ \hline
 \end{array},\;
\begin{array}{|c|c|c|c|c|c|c|c|} \hline 
  -10 & 8 & 2 & 7 & -7 & 19 & -19 \\ \hline  
  6 & 59 & -65 & -91 & 91 & -97 & 97 \\ \hline  
 4 & -68 & 63 & 85 & -84 & 78 & -78 \\ \hline 
 -3 & 89 & -86 & 45 & -44 & 51 & -52 \\ \hline 
 3 & -88 & 86 & -46 & 44 & -51 & 52 \\ \hline
 \end{array},$

 $ \begin{array}{|c|c|c|c|c|c|c|c|} \hline   
  40 & -18 & -22 & 21 & -21 & 23 & -23 \\ \hline   
  -14 & 48 & -34 & -98 & 98 & -99 & 99 \\ \hline 
  -26 & -30 & 56 & 77 & -77 & 76 & -76 \\ \hline
   9 & -92 & 83 & 60 & -62 & -64 & 66 \\ \hline      
      -9 & 92 & -83 & -60 & 62 & 64 & -66 \\ \hline 
 \end{array},\;
 \begin{array}{|c|c|c|c|c|c|c|c|} \hline   
 -40 & 18 & 22 & 25 & -25 & 27 & -27 \\ \hline  
 14 & -48 & 34 & -100 & 100 & -101 & 101 \\ \hline
 26 & 30 & -56 & 75 & -75 & 74 & -74 \\ \hline 
 11 & -93 & 82 & 55 & -57 & -67 & 69 \\ \hline 
 -11 & 93 & -82 & -55 & 57 & 67 & -69 \\ \hline
 \end{array},$

 $\begin{array}{|c|c|c|c|c|c|c|c|} \hline   
     36 & -16 & -20 & 29 & -29 & 31 & -31 \\ \hline   
   -12 & 53 & -41 & -102 & 102 & -103 & 103 \\ \hline   
   -24 & -37 & 61 & 73 & -73 & 72 & -72 \\ \hline   
   13 & -94 & 81 & 39 & -43 & -54 & 58 \\ \hline   
  -13 & 94 & -81 & -39 & 43 & 54 & -58 \\ \hline   
 \end{array},\;
  \begin{array}{|c|c|c|c|c|c|c|c|} \hline   
  -36 & 16 & 20 & 33 & -33 & 35 & -35 \\ \hline   
  12 & -53 & 41 & -104 & 104 & -105 & 105 \\ \hline   
  24 & 37 & -61 & 71 & -71 & 70 & -70 \\ \hline   
  15 & -95 & 80 & 28 & -32 & -38 & 42 \\ \hline   
  -15 & 95 & -80 & -28 & 32 & 38 & -42 \\ \hline   
 \end{array}.$
 \end{footnotesize}
 
\caption{An $\SMAS(5,7;6)$.}\label{576} 
\end{figure}

\begin{lem}\label{57-3}
There exists an $\SMAS(5,7;2c)$ for all $c\geqslant 3$ such that $c\equiv 3 \pmod 4$.
\end{lem}

\begin{proof}
An $\SMAS(5,7;6)$ is shown in Figure \ref{576}. So, write $c=4t+7$ where $t\geqslant 0$. Our $\SMAS(5, 7;2c)$ will consist of $8t+12$ arrays of type \eqref{570} and two arrays of type  
\eqref{Y5}.
Hence, we will construct
$8t+12$ blocks $X_3$, 
$24t+38 $  blocks $X_2$,
$8t+12$   blocks  $Z_4$, $2$ blocks $Y_5$ and
$2$ blocks $W_2$.

Take
$$\begin{array}{rcl}
\overline{\mathfrak{A}} & =& \overline{\mathfrak{A}}( 4t+9,\; 24t+42,\; 10t+17,\; 2t+3 ),\\
\mathfrak{B} & =& \mathfrak{B}(9 ,\; 116t+199 ,\;  24t+38 ),\\
\mathfrak{C} & =& \mathfrak{C}(92t+148,\;116t+200  ,\; 92t+158)
 \end{array}$$
as in Lemmas \ref{BlockA-2}, \ref{BlockB} and \ref{BlockC}, respectively.
We have $$\E(\overline{\mathfrak{A}}) \dcup\E(\mathfrak{B})\dcup \E(\mathfrak{C})=\pm [1, 140t+245]
\setminus \pm \left(\mathcal{M}_1\dcup \ldots \dcup \mathcal{M}_{18} \right),$$
 where
 $$\begin{array}{rclcrcl}
\mathcal{M}_1 & =& [12, 8t+18]_2, & \quad &
\mathcal{M}_2 & =& [40t+68, 40t+72]_4, \\
\mathcal{M}_3 & =& [52t+91,  56t+93]_2, &&
\mathcal{M}_4 & =& [56t+94, 56t+95],\\
\mathcal{M}_5 & =& [56t+98, 60t+100]_2, &&
\mathcal{M}_6 & =& [ 60t+102,  64t+111],\\
\mathcal{M}_7 & =& [64t+113, 64t+115]_2, &&
\mathcal{M}_8 & =& [64t+119, 72t+123]_4,\\
\mathcal{M}_9 & =& [68t+118,  72t+124]_2, &&
\mathcal{M}_{10} & =& [72t+126, 80t+131],\\
\mathcal{M}_{11} & =& [ 80t+132,  80t+135], &&
\mathcal{M}_{12} & =& [80t+143, 88t+147]_4,\\
\mathcal{M}_{13} & =& [80t+137, 80t+139]_2, &&
\mathcal{M}_{14} & =& [84t+142, 92t+144]_2,\\
\mathcal{M}_{15} & =& [88t+151, 92t+149]_2, &&
\mathcal{M}_{16} & =& [92t+146, 92t+150]_4,\\
\mathcal{M}_{17} & =& [92t+151, 92t+153]_2, &&
\mathcal{M}_{18} & =& [92t+155, 92t+156].
\end{array}$$
So, we replace the $8t+12$ instances of $X_3$ with the arrays in  $\overline{\mathfrak{A}}$,  
the $2$  instances of $Y_5$ with the arrays in  $\mathfrak{C}$ and the $24t+38$ instances of $X_2$ with the arrays in $\mathfrak{B}$.

Call
$$\mathfrak{R} =\{U_2(\mathcal{M}_4), U_2(\mathcal{M}_{18}) \}.$$
Then, $\E(\mathfrak{R})=\pm \left(\mathcal{M}_4\dcup \mathcal{M}_{18}\right)$: 
we replace the $2$ instances of $W_2$ with the arrays in $\mathfrak{R}$.

The set $\mathcal{M}_2 \dcup \mathcal{M}_8 \dcup\mathcal{M}_{12}\dcup \mathcal{M}_{16}$  can be written as a disjoint union 
$F_1\dcup \ldots \dcup F_{2\alpha}$, where $\alpha=t+2$ and each $F_i$ is a $2$-set of type $4$.
The set $\mathcal{M}_1\dcup \mathcal{M}_3\dcup \mathcal{M}_5\dcup \mathcal{M}_7 \dcup \mathcal{M}_9\dcup
\mathcal{M}_{13}\dcup \mathcal{M}_{14}\dcup \mathcal{M}_{15}\dcup \mathcal{M}_{17}$  can be written as a disjoint union 
$G_1\dcup \ldots \dcup G_{2\beta}$, where $\beta=4t+5$ and each $G_j$ is a $2$-set of type $2$.
Finally, the set $\mathcal{M}_6\dcup \mathcal{M}_{10}\dcup \mathcal{M}_{11}$
can be written as a disjoint union  $H_1\dcup \ldots \dcup H_{2\gamma}$, 
where $\gamma=3t+5$ and each $H_h$ is a $2$-set of type $1$.
Call
$$\begin{array}{rcl}
\mathfrak{S}_1 & =& \left\{Q_3(F_{2i-1},F_{2i} ): i \in [1,t+2] \right\},\\
\mathfrak{S}_2 & =& \left\{Q_2(G_{2j-1},G_{2j} ): j \in [1,4t+5] \right\},\\
\mathfrak{S}_3 & =& \left\{Q_1(H_{2h-1},H_{2h} ): h \in [1,3t+5] \right\}.
\end{array}$$
  Then $\mathfrak{S}=\mathfrak{S}_1\dcup\mathfrak{S}_2\dcup \mathfrak{S}_3 $ consists of
  $8t+12 $ zero-sum blocks of size $2\times 4$ and
$\E(\mathfrak{S})=\pm  (\mathcal{M}_1\dcup \mathcal{M}_2\dcup \mathcal{M}_3\dcup 
\mathcal{M}_5\dcup \ldots\dcup \mathcal{M}_{17})$.
Hence, we replace the instances of $Z_4$ with the arrays in~$\mathfrak{S}$.
\end{proof}

\begin{lem}\label{BlockC-2}
Given three positive integers $\alpha,\beta,\gamma$ such that 
$\beta>\alpha+15> \gamma+18>40$, there exists 
a set $\overline{\mathfrak{C}}=\overline{\mathfrak{C}}(\alpha,\beta,\gamma)$ consisting of two square blocks of size $5$ such that
 $$\sigma_r(C)=( 0,0,0,1,-1) \equad \sigma_c(C)=(0,0,0,0,0) \text{ for all } C \in \overline{\mathfrak{C}}$$
 and
 $\E(\overline{\mathfrak{C}})$ consists of the following disjoint subsets:
 $$\begin{array}{rcl}
\E(\overline{\mathfrak{C}})  & = & \pm [8,15] \dcup \pm \{22\} \dcup \pm [\gamma, \gamma+3]\dcup \pm \{\alpha,\alpha+4, \alpha+12, \alpha+15\} \dcup \\
&& \pm [\beta, \beta+3]\dcup \pm [\beta+12, \beta+15]. 
       \end{array}$$
\end{lem}

\begin{proof}
Let $C_1$ and $C_2$ be the following arrays:
$$\begin{array}{|c|c|c|c|c|}\hline
 22 &      -14    &          -8    & 13           &    -13 \\ \hline
-12 & -\alpha     &  \alpha+12     & -(\beta+14)  & \beta+14 \\ \hline
-10 & \alpha+15   &   -(\alpha+4)  &  \beta      & -(\beta+1) \\ \hline
 9 & -(\beta+13)  &   \beta+3      & -(\gamma+1) & \gamma+3  \\ \hline
-9 &   \beta+12   & -(\beta+3)     & \gamma+2 &  -(\gamma+3) \\ \hline
  \end{array},$$
$$\begin{array}{|c|c|c|c|c|}\hline
-22 &      14    &          8   &  15      &  -15 \\ \hline
12 & \alpha   & -(\alpha+12)    & -(\beta+15)& \beta+ 15 \\ \hline
10 & -(\alpha+15)  &  \alpha+4  & \beta +1  & -\beta \\ \hline
-11 & \beta+13  &   -(\beta+2)  & \gamma+1 &  -\gamma \\ \hline
 11 & -(\beta+12)  & \beta+2    & -(\gamma+2) & \gamma\\ \hline
  \end{array}.$$
The set $\overline{\mathfrak{C}}=\{C_1,C_2\}$ has the required properties.
  \end{proof}

\begin{figure}[ht]
\begin{footnotesize}
$\begin{array}{|c|c|c|c|c|c|c|c|} \hline   
53  & 3 & -56 & 17 & -17 & 19 & -19 \\\hline
-57 & 59 & -2 & -65 & 65 & -66 & 66 \\\hline 
55  & -60 & 5 & 48 & -48 & 47 & -47 \\\hline
-58 & 4 & 54 & 16 & -18 & -26 & 28 \\\hline
7   & -6 & -1 & -16 & 18 & 26 & -28 \\\hline 
 \end{array},\;
 \begin{array}{|c|c|c|c|c|c|c|c|} \hline  
  -53 & -3 & 56 & 21 & -21 & 23 & -23 \\\hline
  57 & -59 & 2 & -67 & 67 & -68 & 68 \\\hline
  -55 & 60 & -5 & 46 & -46 & 45 & -45 \\\hline 
  58 & -4 & -54 & 31 & -33 & -34 & 36 \\\hline  
  -7 & 6 & 1 & -31 & 33 & 34 & -36 \\\hline
 \end{array},$
 
 $\begin{array}{|c|c|c|c|c|c|c|c|} \hline   
  22 & -14 & -8 & 13 & -13 & 25 & -25 \\\hline
  -12 & -20 & 32 & -63 & 63 & -69 & 69 \\\hline
  -10 & 35 & -24 & 49 & -50 & 44 & -44 \\\hline
  9 & -62 & 52 & -38 & 40 & 29 & -30 \\\hline
 -9 & 61 & -52 & 39 & -40 & -29 & 30 \\\hline
 \end{array},\;
\begin{array}{|c|c|c|c|c|c|c|c|} \hline   
  -22 & 14 & 8 & 15 & -15 & 27 & -27 \\\hline
  12 & 20 & -32 & -64 & 64 & -70 & 70 \\\hline
 10 & -35 & 24 & 50 & -49 & 43 & -43 \\\hline
 -11 & 62 & -51 & 38 & -37 & 41 & -42 \\\hline
 11 & -61 & 51 & -39 & 37 & -41 & 42 \\\hline
\end{array}.$
 \end{footnotesize}
 
\caption{An $\SMAS(5,7;4)$.}\label{574}
\end{figure}

\begin{lem}\label{57-2}
There exists an $\SMAS(5,7;2c)$ for all $c\geqslant 2$ such that $c\equiv 2 \pmod 4$.
\end{lem}

\begin{proof}
An $\SMAS(5,7;4)$ is shown in Figure \ref{574}. So, write $c=4t+6$ where $t\geqslant 0$. 
Our $\SMAS(5,7;2c)$ will consist of $8t+8$ arrays of type \eqref{570}, two arrays of type \eqref{Y5}
and two of type
$$\begin{array}{|c|c|c|c|}\hline
\multicolumn{2}{|c|}{}  &    X_{2}^\T
&     X_{2}^\T\\ \cline{3-4}
 \multicolumn{2}{|c|}{\;\;\smash{\raisebox{.5\normalbaselineskip}{$X_5$}}\;\;\;} & \multicolumn{2}{c|}{Z_4}   \\ \hline  
 \end{array},$$
where each $X_r$ is a zero-sum block of size $r\times 3$ and each $Z_4$ is a zero-sum block of size $2\times 4$.   
Hence, we will construct
$8t+8$ blocks $X_3$, 
$2$ blocks $X_5$, 
$24t+30$  blocks $X_2$,
$8t+10$ blocks  $Z_4$, $2$ blocks $Y_5$ and
$2$ blocks $W_2$.

Take
 $$\begin{array}{rcl}
 \overline{\mathfrak{A}} & =& \overline{\mathfrak{A}}( 4t+12 ,\; 24t+38 ,\; 10t+14,\; 2t+2 ),\\
 \mathfrak{B} & =& \mathfrak{B}(17 ,\; 116t+164,\;  24t+30 ),\\
 \overline{\mathfrak{C}} & =& \overline{\mathfrak{C}}( 92t+118 ,\; 116t+165  ,\;  60t+89  )
  \end{array}$$
 as in Lemmas \ref{BlockA-2}, \ref{BlockB} and \ref{BlockC-2}, respectively.
Furthermore, take $\mathfrak{D}=\{+D,-D\}$, where
$$D=\begin{array}{|c|c|c|}\hline
 116t+169 & 3 & -( 116t+172) \\ \hline
 -( 116t+173) &  116t+175   & -2\\ \hline
  116t+171 & -( 116t+176) & 5\\ \hline
-(116t+174)  & 4 &  116t+170\\ \hline
7 & -6 & -1\\ \hline
     \end{array}.$$    
 Then,  $$\E(\overline{\mathfrak{A}}) \dcup\E(\mathfrak{B})\dcup \E(\overline{\mathfrak{C}})\dcup\E( \mathfrak{D})=\pm [1, 140t+210]
\setminus \pm \left(\mathcal{M}_1\dcup \ldots \dcup \mathcal{M}_{18} \right),$$
where
 $$\begin{array}{rclcrcl}
\mathcal{M}_1 & =& [16,18]_2, & \quad &
\mathcal{M}_2 & =& [20,24]_4,\\
\mathcal{M}_3 & =& [26,8t+24]_2, &&
\mathcal{M}_4 & =& [40t+60, 40t+72]_4,\\
\mathcal{M}_5 & =& [56t+74, 56t+86]_4, &&
\mathcal{M}_6 & =& [52t+81, 56t+83]_2,\\
\mathcal{M}_7 & =& [56t+90, 60t+88]_2, &&
\mathcal{M}_8 & =& [60t+93, 64t+96],\\
\mathcal{M}_9 & =& [64t+97, 64t+107]_2, &&
\mathcal{M}_{10} & =& [64t+109, 72t+113]_4,\\
\mathcal{M}_{11} & =& [68t+102, 80t+112]_2 , &&
\mathcal{M}_{12} & =& [72t+117, 80t+123]_2,\\
\mathcal{M}_{13} & =& [80t+125, 88t+129]_4, &&
\mathcal{M}_{14} & =& [84t+118, 92t+116]_2,\\
\mathcal{M}_{15} & =& [88t+133, 92t+131]_2, &&
\mathcal{M}_{16} & =& [92t+120, 92t+124]_4,\\
\mathcal{M}_{17} & =& [92t+126, 92t+128]_2, &&
\mathcal{M}_{18} & =& [92t+132, 92t+134]_2.
\end{array}$$
So, we replace the $8t+8$ instances of $X_3$ with the arrays in  $\overline{\mathfrak{A}}$,  
the $2$  instances of $Y_5$ with the arrays in  $\overline{\mathfrak{C}}$,
the $2$  instances of $X_5$ with the arrays in  $\mathfrak{D}$ and
the $24t+30$ instances of $X_2$ with the arrays in $\mathfrak{B}$.

Write $\mathcal{M}_8=\mathcal{M}_8^1\dcup \mathcal{M}_8^2\dcup \mathcal{M}_8^3$, where 
$$\mathcal{M}_8^1=[60t+93, 60t+94], \quad \mathcal{M}_8^2=[60t+95, 60t+96] \equad \mathcal{M}_8^3=[60t+97, 64t+96].$$
Define 
$$\mathfrak{R}=\{U_2(\mathcal{M}_{8}^1), U_2(\mathcal{M}_{8}^2) \}.$$
Then, $\E(\mathfrak{R})=\pm(\mathcal{M}_{8}^1\dcup\mathcal{M}_{8}^2)$: we replace the $2$ instances of $W_2$ with the arrays
in $\mathfrak{R}$. 

The set $\mathcal{M}_2 \dcup \mathcal{M}_4 \dcup\mathcal{M}_{5}\dcup \mathcal{M}_{10}\dcup \mathcal{M}_{13}\dcup \mathcal{M}_{16}$ 
can be written as a disjoint union 
$F_1\dcup \ldots \dcup F_{2\alpha}$, where $\alpha=t+4$ and each $F_i$ is a $2$-set of type $4$.
The set $\mathcal{M}_1\dcup \mathcal{M}_3\dcup \mathcal{M}_6\dcup \mathcal{M}_7 \dcup 
\mathcal{M}_9\dcup \mathcal{M}_{11}\dcup \mathcal{M}_{12}\dcup \mathcal{M}_{14}\dcup \mathcal{M}_{15}\dcup \mathcal{M}_{17}\dcup \mathcal{M}_{18}$  can be written as a disjoint union 
$G_1\dcup \ldots \dcup G_{2\beta}$, where $\beta=6t+6$ and each $G_j$ is a $2$-set of type $2$.
Finally, the set $\mathcal{M}_8^3$ 
can be written as a disjoint union  $H_1\dcup \ldots \dcup H_{2t}$, where  each $H_h$ is a $2$-set of type $1$.
Call
$$\begin{array}{rcl}
\mathfrak{S}_1 & =& \left\{Q_3(F_{2i-1},F_{2i} ): i \in [1,t+4] \right\},\\
\mathfrak{S}_2 & =& \left\{Q_2(G_{2j-1},G_{2j} ): j \in [1,6t+6] \right\},\\
\mathfrak{S}_3 & =& \left\{Q_1(H_{2h-1},H_{2h} ): h \in [1,t] \right\}.
 \end{array}$$
Then $\mathfrak{S}=\mathfrak{S}_1\dcup\mathfrak{S}_2\dcup \mathfrak{S}_3 $ consists of
$8t+10$ zero-sum blocks of size $2\times 4$ and
$\E(\mathfrak{S})=\pm  (\mathcal{M}_1\dcup \ldots\dcup\mathcal{M}_7\dcup \mathcal{M}_8^3\dcup 
 \mathcal{M}_9\dcup \ldots\dcup \mathcal{M}_{18})$.
Hence, we replace the instances of $Z_4$ with the arrays in~$\mathfrak{S}$.
\end{proof}

\section{The case \texorpdfstring{$a=6$}{a = 6}}\label{s6}

We now consider the construction of
an $\SMAS(6,b;c)$ when $b\geqslant 5$ is an odd integer.

\begin{lem}\label{sh}
If there exists an $\SMAS(6,b;c)$ for $b,c\geq 1$, then there exists an $\SMAS(6,b+4;c)$.
\end{lem}

\begin{proof}
For every $j\geqslant 0$, let
$$S_j=\begin{array}{|c|}\hline
  Q_1([j+1,j+2],[j+3,j+4])\\\hline
  Q_1([j+5,j+6],[j+7,j+8])\\\hline
  Q_1([j+9,j+10],[j+11,j+12])\\\hline
  \end{array}.$$
Then, $S_j$ is a $6\times 4$ zero-sum block
such that $\E(S_j)=\pm [j+1,j+12]$.
Let $\{R_0,\ldots, R_{c-1}\}$ be an $\SMAS(6,b;c)$.
For every $i \in [0,c-1]$, define the $ 6\times (b+4)$ array
$$T_i=\begin{array}{|c|c|}\hline
R_i & S_{3bc +12i} \\ \hline
\end{array}.$$
Then, the set $\{T_0,\ldots,T_{c-1}\}$ is an $\SMAS(6,b+4;c)$.
\end{proof}

\begin{lem}\label{657-0}
There exist an $\SMAS(6,5;c)$ and an $\SMAS(6,7;c)$  for all  $c\geqslant 4$ such that $c\equiv 0 \pmod 4$.
\end{lem}

\begin{proof}
Write $c=4t+4$, where $t\geqslant 0$.
Our $\SMAS(6,5;c)$ will consist of $4t+4$ blocks of type
\begin{equation}\label{65-0}
\begin{array}{|c|c|}\hline
   X_{2}     &   \\   \cline{1-1}
   X_{2}    &   \\ \cline{1-1}
   X_{2}    & \smash{\raisebox{\normalbaselineskip}{$Z_6^\T$}} \\ \hline  
\end{array},
\end{equation}
whereas our $\SMAS(6,7;c)$ will consist of $4t+4$ blocks of type
\begin{equation}\label{67-0}
\begin{array}{|c|c|c|}\hline
   X_{2}     & \multicolumn{2}{c|}{\;\; Z_4 \;\;}  \\   \hline
   X_{2}    &  \multicolumn{2}{c|}{\;\; Z_4 \;\;} \\ \hline  
   X_{2}    &  \multicolumn{2}{c|}{\;\; Z_4 \;\;} \\ \hline  
\end{array},
\end{equation}
where each $X_2$ is a zero-sum block of size $2\times 3$ and each $Z_\ell$ is a zero-sum block of size $2\times \ell$.
Write  $b=5+\delta $, where $\delta \in \{0, 2\}$, and take
$$\mathfrak{B}=\mathfrak{B}(1,\; 36t+36  , \; 12t+12 )$$
as in Lemma \ref{BlockB}. Then,
$\E(\mathfrak{B})=\pm [1, 12(t+1)(5+\delta)]\setminus \pm (\mathcal{M}_1\dcup \mathcal{M}_2)$,
where
$$\mathcal{M}_1 = [2,24t+24]_2 \equad 
\mathcal{M}_2 = [48t+49, 60t+60+ 12(t+1)\delta ].$$
So, we replace the $12t+12$ instances  of $X_2$ with the arrays in $\mathfrak{B}$.

Write $$\mathcal{M}_1=\mathcal{M}_1^1\dcup\mathcal{M}_1^2\dcup \mathcal{M}_1^3 \equad
\mathcal{M}_2=\mathcal{M}_2^1\dcup\mathcal{M}_2^2\dcup \mathcal{M}_2^3,$$ where
$$\begin{array}{rclcrcl} 
\mathcal{M}_1^1 & =& [2, 8t+6]_4, & \quad & 
\mathcal{M}_2^1 & =& [48t+49, 52t+51]_2,\\
\mathcal{M}_1^2 & =& [4, 8t+8 ]_4, &&
\mathcal{M}_2^2 & =& [48t+50, 52t +52]_2, \\
\mathcal{M}_1^3 & =&  [8t+10,  24t+24]_2, &&
\mathcal{M}_2^3 & =&  [52t+53, 60t+60 + 12(t+1)\delta ].
\end{array}$$
The set $\mathcal{M}_1^1\dcup \mathcal{M}_1^2$ can be written as a disjoint union
$F_1\dcup \ldots \dcup F_\alpha$, where $\alpha=2t+2$ and each $F_i$ is a $2$-set of type $4$.
The set $\mathcal{M}_1^3\dcup \mathcal{M}_2^1\dcup \mathcal{M}_2^2$ can be written as a disjoint union
$G_1\dcup \ldots \dcup G_\beta$, where $\beta=6t+6$ and each $G_j$ is a $2$-set of type $2$.
Finally, the set $\mathcal{M}_2^3$ can be written as a disjoint union
$H_1\dcup \ldots \dcup H_\gamma$, where 
$\gamma=2(2+ 3\delta)(t+1)$ and each $H_h$ is a $2$-set of type $1$.

Suppose $\delta=0$ and call
$$\begin{array}{rcl}
\mathfrak{R}_1 & =& \{R_2(F_i,G_{2i-1},G_{2i} ) :  i \in [1, 2t+2]\},\\
\mathfrak{R}_2 & =& \{R_1(G_{4t+4+j},H_{2j-1},H_{2j} ) :  j \in [1, 2t+2]\}.
\end{array}$$
Then $\mathfrak{R}=\mathfrak{R}_1\dcup \mathfrak{R}_2$ consists of $4t+4$ zero-sum blocks of size $2\times 6$ and 
$\E(\mathfrak{R})=\pm (\mathcal{M}_1\dcup \mathcal{M}_2)$. Hence, we replace the $4t+4$ instances of $Z_6$
with the arrays in $\mathfrak{R}$.

Suppose $\delta = 2$ and call
$$\begin{array}{rcl}
\mathfrak{R}_1 & =& \{Q_3(F_{2i-1},F_{2i} ) :  i \in [1, t+1]\},\\
\mathfrak{R}_2 & =& \{Q_2(G_{2j-1},G_{2j} ) :  j \in [1, 3t+3]\},\\
\mathfrak{R}_3 & =& \{Q_1(H_{2h-1},H_{2h} ) :  h \in [1, 8t+8]\}.
\end{array}$$
Then, $\mathfrak{R}=\mathfrak{R}_1\dcup \mathfrak{R}_2\dcup \mathfrak{R}_3$ consists of $12t+12$ zero-sum blocks of size $2\times 4$ and 
$\E(\mathfrak{R})=\pm (\mathcal{M}_1\dcup \mathcal{M}_2)$. Hence, we replace the $3(4t+4)$ instances of $Z_4$
with the arrays in $\mathfrak{R}$.
\end{proof}

\begin{lem}\label{65-1}
There exists an $\SMAS(6,5;c)$ for all  $c\geqslant 1$ such that $c\equiv 1 \pmod 4$.
\end{lem}

\begin{proof}
Write  $c=4t+1$, where $t\geqslant 0$.
Our $\SMAS(6,5;c)$ will consist of $4t+1$ blocks of type \eqref{65-0}.
Take
$$\mathfrak{B}=\mathfrak{B}(1,\; 36t+8 , \; 12t+3 )$$
as in Lemma \ref{BlockB}. Then $\E(\mathfrak{B})=\pm [1, 15(4t+1) ]\setminus \pm (\mathcal{M}_1\dcup \mathcal{M}_2)$,
where $$\mathcal{M}_1 = [2, 24t+4]_2 \equad \mathcal{M}_2 = [ 48t+12,  60t+15 ].$$
Write $$\mathcal{M}_1=\mathcal{M}_1^1\dcup\mathcal{M}_1^2\dcup \mathcal{M}_1^3 \equad
\mathcal{M}_2=\mathcal{M}_2^1\dcup\mathcal{M}_2^2\dcup \mathcal{M}_2^3,$$ where
$$\begin{array}{rclcrcl}
\mathcal{M}_1^1 & =& [2, 16t+4]_2, & \quad &
\mathcal{M}_2^1 & =& [48t+12, 52t+10]_2,\\
\mathcal{M}_1^2 & =& [16t+6,   24t+2 ]_4, &&
\mathcal{M}_2^2 & =& [48t+13, 52t +11]_2, \\
\mathcal{M}_1^3 & =&  [16t+8,   24t+4]_{4}, &&
\mathcal{M}_2^3 & =&  [52t+12, 60t+15].
\end{array}$$
So, we replace the $12t+3$ instances  of $X_2$ with the arrays in $\mathfrak{B}$.

The set $\mathcal{M}_1^2\dcup \mathcal{M}_1^3$ can be written as a disjoint union
$F_1\dcup \ldots \dcup F_\alpha$, where $\alpha=2t$ and each $F_i$ is a $2$-set of type $4$.
The set $\mathcal{M}_1^1\dcup \mathcal{M}_2^1\dcup \mathcal{M}_2^2$ can be written as a disjoint union
$G_1\dcup \ldots \dcup G_\beta$, where $\beta=6t+1 $ and each $G_j$ is a $2$-set of type $2$.
Finally, the set $\mathcal{M}_2^3$ can be written as a disjoint union
$H_1\dcup \ldots \dcup H_\gamma$, where 
$\gamma=4t+2$ and each $H_h$ is a $2$-set of type $1$.
Call
$$\begin{array}{rcl}
\mathfrak{R}_1 & =& \{R_2(F_i,G_{2i-1},G_{2i} ) :  i \in [1, 2t]\},\\
\mathfrak{R}_2 & =& \{R_1(G_{4t+j},H_{2j-1},H_{2j} ) :  j \in [1, 2t+1]\}.
\end{array}$$
Then $\mathfrak{R}=\mathfrak{R}_1\dcup \mathfrak{R}_2$ consists of $4t+1$ zero-sum blocks of size $2\times 6$ and 
$\E(\mathfrak{R})=\pm (\mathcal{M}_1\dcup \mathcal{M}_2)$. Hence, we replace the $4t+1$ instances of $Z_6$
with the arrays in $\mathfrak{R}$.
\end{proof}

\begin{lem}\label{67-1}
There exists an $\SMAS(6,7;c)$ for all  $c\geqslant 1$ such that $c\equiv 1 \pmod 4$.
\end{lem}

\begin{proof}
The existence of an $\SMAS(6,7;1)$ follows from \cite[Lemma 7]{KCW}.
So, write  $c=4t+5$, where $t\geqslant 0$.
Our $\SMAS(6,7;c)$ will consist of $4t+4$ blocks of type \eqref{67-0}
and one block
\begin{equation}\label{X6}
\begin{array}{|c|c|c|}\hline
 &   \multicolumn{2}{c|}{Z_4}\\   \cline{2-3}
  &  \multicolumn{2}{c|}{Z_4} \\ \cline{2-3}
\smash{\raisebox{\normalbaselineskip}{$X_6$}}  & \multicolumn{2}{c|}{Z_4}\\\hline
\end{array},
\end{equation}
where $X_6$ is a zero-sum block of size $6\times 3$ and each $Z_4$ is a zero-sum block of size $2\times 4$.
Take the zero-sum block
\begin{equation}\label{specA}
A=\begin{array}{|c|c|c|c|c|c|}\hline
1 & -1 & 2 & -2 & 4 & -4 \\\hline
8 & 7 & -8 & 5 & -7 & -5 \\\hline
-9 & -6 & 6 & -3 & 3 & 9 \\\hline
\end{array} 
\end{equation}
and the set 
$$\mathfrak{B}=\mathfrak{B}(11,\;  36t+45, \; 12t+12 )$$
as in Lemma \ref{BlockB}. Then $\E(A)\dcup \E(\mathfrak{B})=\pm [1, 21(4t+5) ]\setminus \pm (\mathcal{M}_1
\dcup \mathcal{M}_2\dcup \mathcal{M}_3)$,
where
$$\begin{array}{rclcrcl}
\mathcal{M}_1 & = & [10, 24t+32]_2 , & \quad &
\mathcal{M}_2  & =&  [36t+46, 36t+55], \\
\mathcal{M}_3 & = & [48t+68, 84t+105].
  \end{array}$$
So, we replace the unique instance of $X_6$ with $A^\T$ and the $12t+12$ instances  of $X_2$ with the arrays in $\mathfrak{B}$.

The set $\mathcal{M}_1$ can be written as a disjoint union
$G_1\dcup \ldots \dcup G_{2\beta}$, where $\beta=3t+3$ and each $G_j$ is a $2$-set of type $2$.
The set $\mathcal{M}_2\dcup \mathcal{M}_3$ can be written as a disjoint union
$H_1\dcup \ldots \dcup H_{2\gamma}$, where $\gamma=9t+12$ and each $H_h$ is a $2$-set of type $1$.
Call
$$\begin{array}{rcl}
\mathfrak{R}_1 & =& \{Q_2(G_{2j-1},G_{2j} ) :  j \in [1, 3t+3]\},\\
\mathfrak{R}_2 & =& \{Q_1(H_{2h-1},H_{2h} ) :  h \in [1, 9t+12]\}.
\end{array}$$
Then, $\mathfrak{R}=\mathfrak{R}_1\dcup \mathfrak{R}_2$ consists of $12t+15$ zero-sum blocks of size $2\times 4$ and 
$\E(\mathfrak{R})=\pm (\mathcal{M}_1\dcup \mathcal{M}_2\dcup \mathcal{M}_3)$. Hence, we replace the $3(4t+5)$ instances of $Z_4$
with the arrays in $\mathfrak{R}$.
\end{proof}

\begin{figure}[ht]
\begin{footnotesize}
$\begin{array}{|c|c|c|c|c|} \hline   
-4 & -5 & 9 & 11 & -11 \\\hline
4 & -7 & 3 & -13 & 13  \\\hline
-2 & 5 & -3 & -24 & 24  \\\hline
2 & -8 & 6 & 25 & -25  \\\hline
-1 & 7 & -6 & -26 & 26  \\\hline
1 & 8 & -9 & 27 & -27  \\\hline
 \end{array},\;
 \begin{array}{|c|c|c|c|c|} \hline   
 10 & -28 & 17 & 19 & -18 \\\hline
-10 & 28 & -17 & -23 & 22 \\\hline
12 & -29 & 16 & -20 & 21 \\\hline
-12 & 29 & -16 & -19 & 18 \\\hline
14 & -30 & 15 & 23 & -22 \\\hline
-14 & 30 & -15 & 20 & -21 \\\hline
 \end{array}$.
 \end{footnotesize}
 
\caption{An $\SMAS(6,5;2)$.}\label{652}
\end{figure}

\begin{lem}\label{65-2}
There exists an $\SMAS(6,5;c)$ for all  $c\geqslant 2$ such that $c\equiv 2 \pmod 4$.
\end{lem}

\begin{proof}
An $\SMAS(6,5;2)$ is shown in Figure \ref{652}.
So, write  $c=4t+6$, where $t\geqslant 0$.
Our $\SMAS(6,5;c)$ will consist of $4t+5$ blocks of type \eqref{65-0}
and one block
\begin{equation}\label{X6-1}
\begin{array}{|c|c|c|}\hline
 &   \\   
  &   \\ 
\smash{\raisebox{\normalbaselineskip}{$X_6$}}  & \smash{\raisebox{\normalbaselineskip}{$Z_6^\T$}}\\\hline
\end{array},
\end{equation}
where $X_6$ is a zero-sum block of size $6\times 3$ and  $Z_6$ is a zero-sum block of size $2\times 6$.
Take the zero-sum block $A$ given in \eqref{specA} 
and the set 
$$\mathfrak{B}=\mathfrak{B}(11,\;  36t+55 , \; 12t+15  )$$
as in Lemma \ref{BlockB}. Then $\E(A)\dcup \E(\mathfrak{B})=\pm [1, 15(4t+6) ]\setminus \pm (\mathcal{M}_1
\dcup \mathcal{M}_2\dcup \mathcal{M}_3)$,
where
$$\begin{array}{rclcrcl}
\mathcal{M}_1 & = & [10, 24t+40]_2 , & \quad &
\mathcal{M}_2  & =&  [36t+56, 36t+65],\\
\mathcal{M}_3 & = &  [48t+81, 60t+90].
  \end{array}$$
So, we replace the unique instance of $X_6$ with $A^\T$ and the $12t+15$ instances  of $X_2$ with the arrays in $\mathfrak{B}$.

Write $$\mathcal{M}_1=\mathcal{M}_1^1\dcup\mathcal{M}_1^2\dcup \mathcal{M}_1^3, \quad
\mathcal{M}_2=\mathcal{M}_2^1\dcup\mathcal{M}_2^2\dcup \mathcal{M}_2^3 \equad
\mathcal{M}_3=\mathcal{M}_3^1\dcup\mathcal{M}_3^2\dcup \mathcal{M}_3^3,$$ where
$$\begin{array}{rclcrcl}
\mathcal{M}_1^1 & =& [10, 8t+24]_2, &  \quad & 
\mathcal{M}_2^1 & =& [36t+56, 36t+58]_2,\\
\mathcal{M}_1^2 & =& [8t+26,  24t+38 ]_4, &&
\mathcal{M}_2^2 & =& [36t+57, 36t+59  ]_2, \\
\mathcal{M}_1^3 & =& [8t+28,  24t+40]_4, &&
\mathcal{M}_2^3 & =&  [36t+60, 36t+65],\\
\mathcal{M}_3^1 & =& [48t+81, 60t+87]_2, &&
\mathcal{M}_3^2 & =& [48t+82, 60t+88  ]_2,\\
\mathcal{M}_3^3 & =&  [60t+89,  60t+90].
\end{array}$$
The set $\mathcal{M}_1^2\dcup \mathcal{M}_1^3$ can be written as a disjoint union
$F_1\dcup \ldots \dcup F_\alpha$, where $\alpha=4t+4$ and each $F_i$ is a $2$-set of type $4$.
The set $\mathcal{M}_1^1\dcup \mathcal{M}_2^1\dcup \mathcal{M}_2^2 \dcup \mathcal{M}_3^1\dcup \mathcal{M}_3^2$ can be written as a disjoint union
$G_1\dcup \ldots \dcup G_{\beta}$, where $\beta=8t+10$ and each $G_j$ is a $2$-set of type $2$.
Finally, the set $\mathcal{M}_2^3\dcup \mathcal{M}_3^3$ can be written as  a disjoint union
$H_1\dcup \ldots \dcup H_{4}$, where each $H_h$ is a $2$-set of type $1$.
Call
$$\begin{array}{rcl}
\mathfrak{R}_1 & =& \{R_2(F_i,G_{2i-1},G_{2i} ) :  i \in [1, 4t+4]\},\\
\mathfrak{R}_2 & =& \{R_1(G_{8t+8+j}, H_{2j-1},H_{2j} ) :  j \in [1, 2]\}.
\end{array}$$
Then, $\mathfrak{R}=\mathfrak{R}_1\dcup \mathfrak{R}_2$ consists of $4t+6$ zero-sum blocks of size $2\times 6$ and 
$\E(\mathfrak{R})=\pm (\mathcal{M}_1\dcup \mathcal{M}_2\dcup \mathcal{M}_3)$. Hence, we replace the $4t+6$ instances of $Z_6$
with the arrays in $\mathfrak{R}$.
\end{proof}

\begin{lem}\label{67-2}
There exists an $\SMAS(6,7;c)$ for all  $c\geqslant 2$ such that $c\equiv  2 \pmod 4$.
\end{lem}

\begin{proof}
Write  $c=4t+2$, where $t\geqslant 0$.
Our $\SMAS(6,7;c)$ will consist of $4t+1$ blocks of type \eqref{67-0}
and one block of type \eqref{X6}.
Take the zero-sum block $A$ given in \eqref{specA} 
and the set 
$$\mathfrak{B}=\mathfrak{B}(11,\; 36t+19  , \; 12t+3 )$$
as in Lemma \ref{BlockB}. Then $\E(A)\dcup \E(\mathfrak{B})=\pm [1, 21(4t+2) ]\setminus \pm (\mathcal{M}_1
\dcup \mathcal{M}_2\dcup \mathcal{M}_3)$,
where
$$\begin{array}{rclcrcl}
\mathcal{M}_1 & = &  [10, 24t+16]_2, & \quad &
\mathcal{M}_2  & =&  [36t+20, 36t+29],\\
\mathcal{M}_3 & = &  [48t+33, 84t+42].
  \end{array}$$
So, we replace the unique instance of $X_6$ with $A^\T$ and the $12t+3$ instances  of $X_2$ with the arrays in $\mathfrak{B}$.

The set $\mathcal{M}_1$ can be written as a disjoint union
$G_1\dcup \ldots \dcup G_{2\beta}$, where $\beta=3t+1$ and each $G_j$ is a $2$-set of type $2$.
The set $\mathcal{M}_2\dcup \mathcal{M}_3$ can be written as a disjoint union
$H_1\dcup \ldots \dcup H_{2\gamma}$, where $\gamma=9t+5$ and each $H_h$ is a $2$-set of type $1$.
Call
$$\begin{array}{rcl}
\mathfrak{R}_1 & =& \{Q_2(G_{2j-1},G_{2j} ) :  j \in [1, 3t+1]\},\\
\mathfrak{R}_2 & =& \{Q_1(H_{2h-1},H_{2h} ) :  h \in [1, 9t+5]\}.
\end{array}$$
Then, $\mathfrak{R}=\mathfrak{R}_1\dcup \mathfrak{R}_2$ consists of $12t+6$ zero-sum blocks of size $2\times 4$ and 
$\E(\mathfrak{R})=\pm (\mathcal{M}_1\dcup \mathcal{M}_2\dcup \mathcal{M}_3)$. Hence, we replace the $3(4t+2)$ instances of $Z_4$
with the arrays in $\mathfrak{R}$.
\end{proof}

\begin{lem}\label{65-3}
There exists an $\SMAS(6,5;c)$ for all  $c\geqslant 3$ such that $c\equiv 3 \pmod 4$.
\end{lem}

\begin{proof}
Write  $c=4t+3$, where $t\geqslant 0$.
Our $\SMAS(6,5;c)$ will consist of $4t+2$ blocks of type \eqref{65-0}
and one block \eqref{X6-1}.
Take the zero-sum block $A$ given in \eqref{specA} 
and the set 
$$\mathfrak{B}=\mathfrak{B}(11,\;  36t+27 , \; 12t+6 )$$
as in Lemma \ref{BlockB}. Then $\E(A)\dcup \E(\mathfrak{B})=\pm [1, 15(4t+3) ]\setminus \pm (\mathcal{M}_1
\dcup \mathcal{M}_2\dcup \mathcal{M}_3)$,
where  
$$\begin{array}{rclcrcl}
\mathcal{M}_1 & = &  [10, 24t+20]_2 , & \quad &
\mathcal{M}_2  & =&  [36t+28, 36t+37],\\
\mathcal{M}_3 & = &  [48t+44, 60t+45].
  \end{array}$$
So, we replace the unique instance of $X_6$ with $A^\T$ and the $12t+6$ instances  of $X_2$ with the arrays in $\mathfrak{B}$.

Write $$\mathcal{M}_1=\mathcal{M}_1^1\dcup\mathcal{M}_1^2\dcup \mathcal{M}_1^3\equad
\mathcal{M}_3=\mathcal{M}_3^1\dcup\mathcal{M}_3^2\dcup \mathcal{M}_3^3,$$ where 
$$\begin{array}{rclcrcl}
\mathcal{M}_1^1 & =& [10, 8t+20]_2, & \quad &
\mathcal{M}_3^1 & =& [48t+44, 60t+42]_2,\\
\mathcal{M}_1^2 & =& [8t+22,  24t+18 ]_4, &&
\mathcal{M}_3^2 & =& [48t+45, 60t+43]_2, \\
\mathcal{M}_1^3 & =& [8t+24,  24t+20]_4, &&
\mathcal{M}_3^3 & =&  [60t+44, 60t+45].
\end{array}$$
The set $\mathcal{M}_1^2\dcup \mathcal{M}_1^3$ can be written as a disjoint union
$F_1\dcup \ldots \dcup F_\alpha$, where $\alpha=4t$ and each $F_i$ is a $2$-set of type $4$.
The set $\mathcal{M}_1^1\dcup \mathcal{M}_3^1\dcup \mathcal{M}_3^2$ can be written as a disjoint union
$G_1\dcup \ldots \dcup G_{\beta}$, where $\beta=8t+3$ and each $G_j$ is a $2$-set of type $2$.
Finally, the set $\mathcal{M}_2\dcup \mathcal{M}_3^3$ can be written as a disjoint union
$H_1\dcup \ldots \dcup H_{6}$, where each $H_h$ is a $2$-set of type $1$.
Call
$$\begin{array}{rcl}
\mathfrak{R}_1 & =& \{R_2(F_i,G_{2i-1},G_{2i} ) :  i \in [1, 4t]\},\\
\mathfrak{R}_2 & =& \{R_1(G_{8t+j}, H_{2j-1},H_{2j} ) :  j \in [1, 3]\}.
\end{array}$$
Then, $\mathfrak{R}=\mathfrak{R}_1\dcup \mathfrak{R}_2$ consists of $4t+3$ zero-sum blocks of size $2\times 6$ and 
$\E(\mathfrak{R})=\pm (\mathcal{M}_1\dcup \mathcal{M}_2\dcup \mathcal{M}_3)$. Hence, we replace the $4t+3$ instances of $Z_6$
with the arrays in $\mathfrak{R}$.
\end{proof}

\begin{lem}\label{67-3}
There exists an $\SMAS(6,7;c)$ for all  $c\geqslant 3$ such that $c\equiv 3 \pmod 4$.
\end{lem}

\begin{proof}
Write  $c=4t+3$, where $t\geqslant 0$.
Our $\SMAS(6,7;c)$ will consist of $4t+3$ blocks of type \eqref{67-0}.
Take
$$\mathfrak{B}=\mathfrak{B}(1,\; 36t+26 , \;  12t+9 )$$
as in Lemma \ref{BlockB}. Then $\E(\mathfrak{B})=\pm [1, 21(4t+3) ]\setminus \pm (\mathcal{M}_1\dcup \mathcal{M}_2)$,
where $$\mathcal{M}_1 = [2, 24t+16]_2 \equad \mathcal{M}_2 = [ 48t+36,  84t+63 ].$$
So, we replace the $12t+9$ instances  of $X_2$ with the arrays in $\mathfrak{B}$.

The set $\mathcal{M}_1$ can be written as a disjoint union
$G_1\dcup \ldots \dcup G_{2\beta}$, where $\beta=3t+2$ and each $G_j$ is a $2$-set of type $2$.
The set $\mathcal{M}_2$ can be written as a disjoint union
$H_1\dcup \ldots \dcup H_{2\gamma}$, where 
$\gamma=9t+7$ and each $H_h$ is a $2$-set of type $1$.
Call
$$\begin{array}{rcl}
\mathfrak{R}_1 & =& \{Q_2(G_{2j-1},G_{2j} ) :  j \in [1, 3t+2]\},\\
\mathfrak{R}_2 & =& \{Q_1(H_{2h-1},H_{2h} ) :  h \in [1, 9t+7]\}.
\end{array}$$
Then, $\mathfrak{R}=\mathfrak{R}_1\dcup \mathfrak{R}_2$ consists of $12t+9$ zero-sum blocks of size $2\times 4$ and 
$\E(\mathfrak{R})=\pm (\mathcal{M}_1\dcup \mathcal{M}_2)$. Hence, we replace the $3(4t+3)$ instances of $Z_4$
with the arrays in $\mathfrak{R}$.
\end{proof}

\section{Conclusions}\label{concl}

We can now prove our main results.

\begin{proof}[Proof of Proposition \ref{p:sets}]
By Remark \ref{tr}, the statement follows from the results of Sections \ref{gen} and \ref{sm}, according to 
Table \ref{Tab}.
\end{proof}

\begin{table}[ht]
\begin{tabular}{|c|c|c|c|}\hline
$\bm{c\equiv 0\pmod 4}$           & $b=7$              & $9\leqslant b \equiv 1 \pmod 4$ & $11\leqslant b \equiv 3 \pmod 4$ \\\hline
$a=5$                             & L. \ref{57-0}      & L. \ref{57-0}                   & L. \ref{57-0}  \\\hline
$a\geqslant 7$                    & C. \ref{IHS->SMAS} & C. \ref{IHS->SMAS}              & C. \ref{IHS->SMAS} \\\hline \hline

$\bm{c\equiv 1\pmod 4}$                & $b=7$              & $9\leqslant b \equiv 1 \pmod 4$ & $11\leqslant b \equiv 3 \pmod 4$ \\\hline
$a=5$                             & L. \ref{57-1}      & L. \ref{a1b1-1}                 & L. \ref{511-1}  \\\hline
$7\leqslant a \equiv 3 \pmod 4$   & L. \ref{a3b3-2}    & C. \ref{IHS->SMAS}              & L. \ref{a3b3-2}  \\\hline
$9\leqslant a \equiv 1 \pmod 4$   & C. \ref{IHS->SMAS} & L. \ref{a1b1-1}                 & C. \ref{IHS->SMAS}  \\\hline\hline

$\bm{c\equiv 2\pmod 4}$                & $b=7$              & $9\leqslant b \equiv 1 \pmod 4$ & $11\leqslant b \equiv 3 \pmod 4$ \\\hline
$a=5$                             & L. \ref{57-2}      & L. \ref{59-2n}                  & L. \ref{59-2n}   \\\hline
$a=7$                             & L. \ref{ab3-2}     & L. \ref{97-c2}                  & L. \ref{ab3-2}              \\\hline
$9\leqslant a \equiv 1 \pmod 4$   & L. \ref{97-c2}     & L. \ref{59-2n}                  & L. \ref{59-2n}             \\\hline
$11\leqslant a \equiv 3 \pmod 4$  & L. \ref{ab3-2}     & L. \ref{59-2n}                  & L. \ref{ab3-2}              \\\hline\hline

$\bm{c\equiv 3\pmod 4}$                & $b=7$              & $9\leqslant b \equiv 1 \pmod 4$ & $11\leqslant b \equiv 3 \pmod 4$ \\\hline
$a=5$                             & L. \ref{57-3}      & L. \ref{59-3}                   & L. \ref{511-3}   \\\hline
$a=7$                             & C. \ref{IHS->SMAS} & L. \ref{97-3}                   & C. \ref{IHS->SMAS} \\\hline
$9\leqslant a \equiv 1 \pmod 4$   & L. \ref{97-3}      & C. \ref{IHS->SMAS}              & L. \ref{511-3}   \\\hline
$11\leqslant a \equiv 3 \pmod 4$  & C. \ref{IHS->SMAS} & L. \ref{511-3}                  & C. \ref{IHS->SMAS} \\\hline
\end{tabular}
\caption{Proof of Proposition \ref{p:sets}.}\label{Tab}
\end{table}

\begin{proof}[Proof of Proposition \ref{6}]
For all $c\geqslant 1$, there exists
an $\SMAS(6,5;c)$ by 
Lemmas \ref{657-0}, \ref{65-1}, \ref{65-2} and \ref{65-3},
and there  exists
an $\SMAS(6,7;c)$ by 
Lemmas \ref{657-0}, \ref{67-1}, \ref{67-2} and \ref{67-3},
Hence, the result follows from Lemma~\ref{sh}.
\end{proof}

\begin{proof}[Proof of Theorem \ref{main}]
By \cite[Main Theorem]{KE} and \cite[Main Theorem 20]{KLE}, we may assume $s,k\geq 4$.
Up to transposition, by Proposition \ref{fatto} we may also assume that $k\geqslant 5$ is odd, 
$n$ is even, $s\not \equiv 0 \pmod 4$ and $\gcd(s,k)=1$.
From $ms=nk$ we get that there exists a positive integer $e$ such that $m=ek$ and $n=es$.

Suppose that $s$ is odd. Since $n$ is even, $e=2c$ for a suitable positive integer $c$.
Let $\{R_1,\ldots,R_{2c}\}$ be an $\SMAS(k,s;2c)$: the existence of this set follows from Proposition \ref{p:sets},
since $s\geqslant 5$ and $(s,k)\neq (5,5)$.
In this case, the $m\times n$ block-diagonal matrix $\diag(R_1,\ldots,R_{2c})$ is an $\SMA(m,n;s,k)$.

Next, suppose that $s$ is even. By the previous assumptions, we can write $s=2\overline{s}$, where $\overline{s}\geqslant 3$
is an odd integer, coprime with $k$.
If $\overline{s}= 3$, let
$\{R_1,\ldots,R_{e}\}$ be an $\SMAS(k,6;e)$, whose existence follows from Proposition \ref{6}.
In this case, the $m\times n$ block-diagonal matrix $\diag(R_1,\ldots,R_{e})$ is an $\SMA(m,n; s,k)$.
If $\overline{s}\geqslant 5$, let $\{R_1,\ldots,R_{2e}\}$ be an $\SMAS(k,\overline{s};2e)$, whose existence  follows from Proposition \ref{p:sets}
since  $(k,\overline{s})\neq (5,5)$.
For every $i\in [1,e]$, define the $k\times s$ array $S_{i} =\begin{array}{|c|c|}\hline R_{2i-1}& R_{2i} \\ \hline \end{array}$.
In this case, the $m\times n$ block-diagonal matrix $\diag(S_1,\ldots,S_{e})$ is an $\SMA(m,n;s,k)$.
\end{proof}

\begin{proof}[Proof of Corollary \ref{fold}]
By \cite[Proposition 4.2]{CP}, there is no (integer) ${}^2\H(m,n;s,k)$ when $nk\equiv 1 \pmod 4$. 
So, assume that $nk$ even.
By Theorem \ref{main} there exists an $\SMA(m,n;s,k)$ (and hence, an integer ${}^2\H(m,n;s,k)$)
when 
$(m,n,s,k)$ does not belong to the set
$$\left\{ (\ell, 2,2,\ell), (2,\ell,\ell, 2): \ell\geqslant 1 \text{ is such that } \ell \equiv 1,2 \pmod 4\right\}.$$
So, we are left to consider the case when
$m=k=2$ and $n=s\equiv 1,2 \pmod 4$. 
For the sake of contradiction, suppose that there exists an integer ${}^2\H(2,n;n,2)$, say $A$.
Denoting by $R_1,R_2$ the two rows of $A$, we must have 
$R_2=-R_1$, where the absolute values of the entries of $R_1$ constitute the set $[1,n]$.
Since $n\equiv 1,2\pmod 4$, the sum of the entries of $R_1$ cannot be equal to zero, a contradiction. 
\end{proof}

\section*{Acknowledgments}

The two authors are members of INdAM-GNSAGA.

\end{document}